\documentclass[11pt, a4paper]{article}
\usepackage[latin1]{inputenc}
\usepackage[T1]{fontenc}
\usepackage[english]{babel}
\usepackage{amssymb}
\usepackage{amsmath}
\usepackage{amsthm}
\usepackage[colorlinks=true, allcolors=blue]{hyperref}
\usepackage{geometry}
\usepackage{fancyhdr}
\usepackage{relsize}
\usepackage{upgreek}
\usepackage{verbatim}
\usepackage{enumitem}
\usepackage{calrsfs}
\usepackage{tikz}
\usetikzlibrary{matrix,arrows,calc}
\usetikzlibrary{decorations.pathreplacing}
\usepackage{mypackage}
\renewcommand{\Hilb}{H}
\newcommand{\treeA}{\text{\hspace{-.7em}
	\begin{tikzpicture}
		\node (s) {};%
		\node (s1) at ($(s) + (-.45em,-.25em)$) {};%
		\node (s11) at ($(s1) + (-.25em,-.25em)$) {};%
		\draw (s.center) -- (s1.center);%
		\draw (s1.center) -- (s11.center);%
		\filldraw (s) circle (.5pt);%
		\filldraw (s1) circle (.5pt);%
		\filldraw (s11) circle (.5pt);%
	\end{tikzpicture}%
}}
\newcommand{\treeB}{\text{\hspace{-.5em}
	\begin{tikzpicture}
		\node (s) {};%
		\node (s1) at ($(s) + (-.45em,-.25em)$) {};%
		\node (s12) at ($(s1) + (.25em,-.25em)$) {};%
		\draw (s.center) -- (s1.center);%
		\draw (s1.center) -- (s12.center);%
		\filldraw (s) circle (.5pt);%
		\filldraw (s1) circle (.5pt);%
		\filldraw (s12) circle (.5pt);%
	\end{tikzpicture}%
}}
\newcommand{\treeC}{\text{\hspace{-.7em}
	\begin{tikzpicture}
		\node (s) {};%
		\node (s1) at ($(s) + (-.4em,-.25em)$) {};%
		\node (s2) at ($(s) + (.4em,-.25em)$) {};%
		\draw (s.center) -- (s1.center);%
		\draw (s.center) -- (s2.center);%
		\filldraw (s) circle (.5pt);%
		\filldraw (s1) circle (.5pt);%
		\filldraw (s2) circle (.5pt);%
	\end{tikzpicture}%
}}
\newcommand{\treeD}{\text{\hspace{-.4em}
	\begin{tikzpicture}
		\node (s) {};%
		\node (s2) at ($(s) + (.45em,-.25em)$) {};%
		\node (s21) at ($(s2) + (-.25em,-.25em)$) {};%
		\draw (s.center) -- (s2.center);%
		\draw (s2.center) -- (s21.center);%
		\filldraw (s) circle (.5pt);%
		\filldraw (s2) circle (.5pt);%
		\filldraw (s21) circle (.5pt);%
	\end{tikzpicture}%
}}
\newcommand{\treeE}{\text{\hspace{-.6em}
	\begin{tikzpicture}
		\node (s) {};%
		\node (s2) at ($(s) + (.45em,-.25em)$) {};%
		\node (s22) at ($(s2) + (.25em,-.25em)$) {};%
		\draw (s.center) -- (s2.center);%
		\draw (s2.center) -- (s22.center);%
		\filldraw (s) circle (.5pt);%
		\filldraw (s2) circle (.5pt);%
		\filldraw (s22) circle (.5pt);%
	\end{tikzpicture}%
}}
\title{On Cohomology Rings of Non-Commutative Hilbert Schemes and CoHa-Modules}
\author{H. Franzen\footnote{\myaddress}}
\date{}
\begin{document}
	\maketitle
	\begin{abstract}
	We prove that Chow groups of certain non-commutative Hilbert schemes have a basis consisting of monomials in Chern classes of the universal bundle. Furthermore, we realize the cohomology of non-commutative Hilbert schemes as a module over the Cohomological Hall algebra.
\end{abstract}

	\section*{Introduction}

The cohomology of the Hilbert scheme of $d$ points in an $m$-dimensional affine space has been studied intensively by various authors (e.g.\ \cite{Grojnowski:96}, \cite{Nakajima:97} and \cite{LS:01}). The objective of this paper is to investigate cohomological properties of a certain non-commutative analog of these Hilbert schemes. Observing that the Hilbert scheme of $d$ points in $\A^m$ parametrizes ideals of codimension $d$ of the polynomial algebra in $m$ variables, we might ask for the moduli space of \emph{left}-ideals of codimension $d$ in the free \emph{non-commutative} algebra in $m$ letters. This is the most prominent example of a non-commutative Hilbert scheme.

So far, we know that non-commutative Hilbert schemes possess a cell decomposition. This was shown by Reineke \cite[Thm.\ 1.3]{Reineke:05}. The cells are parametrized by $m$-ary trees with $d$ nodes. The existence of a cell decomposition implies that the Chow group (and also the singular cohomology) is a free group with a basis given by the closures of the cells. As a consequence, it is possible to give an explicit formula for the Poincar\'{e} polynomial and for the Euler characteristic. But as a non-commutative Hilbert scheme is also a non-singular variety, we know that its Chow group possesses a ring structure. It turns out that calculating the intersection product of the cell closures is a difficult task and therefore, this basis is not so well-adapted to the multiplication. We provide another basis of the Chow ring that allows us insights into the multiplicative structure.

Again, let's turn our attention to classical Hilbert schemes for a moment. A result of Lehn and Sorger (cf.\ \cite[Thm.\ 1.1]{LS:01}) shows that the cohomology ring of the Hilbert scheme of $d$ points in the affine plane is isomorphic to the ring of class functions of the symmetric group $S_d$. This is done by using results of Grojnowski (cf.\ \cite{Grojnowski:96}) and Nakajima (cf.\ \cite{Nakajima:97}) which show that the direct sum (over all $d$) of all these cohomology groups (tensored to the rationals) has the structure of a vertex algebra isomorphic to the bosonic Fock space. It turns out that an appropriate analog in the non-commutative case would be to give the cohomology of non-commu\-ta\-tive Hilbert schemes a module structure over Kontsevich--Soibelman's Cohomological Hall algebra (cf.\ \cite{KS:11}).

The two major results of this paper are the following: As a non-commutative Hilbert scheme arises as a fine moduli space, it is equipped with a universal bundle. We will exhibit a basis of the Chow group consisting of monomials in the Chern classes of the universal bundle in Theorem~\ref{ko_basis_Chern}. In particular, this gives a description of the Chow ring as a quotient of a polynomial ring. We thus circumvent the problem that a result of King and Walter \cite[Thm.\ 3]{KW:95} does not apply here. Their theorem states that Chow rings of fine quiver moduli are generated by Chern classes of tautological bundles if the quiver is \emph{acyclic}. Moreover, we realize the cohomology of non-commutative Hilbert schemes (which equals their Chow ring after extending scalars to the rationals) as a quotient of the Cohomological Hall algebra and describe the kernel of the quotient map explicitly. This is Theorem \ref{thm_kernel}.

The paper is organized as follows: In the first section, we recollect some facts on non-commutative Hilbert schemes that are essential for our purposes. In particular, we present Reineke's cell decomposition (cf.\ Theorem~\ref{thm_cell}). The first main result of this paper, the existence of a basis of the Chow ring of the non-commutative Hilbert scheme consisting of monomials in Chern classes of the universal bundle, is Theorem \ref{ko_basis_Chern}. It is a direct consequence of Theorem \ref{thm} which states that certain monomials in Chern classes, parametrized by $m$-ary trees (or $m$-ary forests for a slightly more general type of non-commutative Hilbert schemes) can be expressed as integer linear combinations of cell closures. These linear combinations provide an upper unitriangular base change matrix. This theorem is proved by expressing the filtration steps of Reineke's cell decomposition as intersections of degeneracy loci and proving that every irreducible component of this intersection has the ``correct'' dimension. The third section is devoted to the description of the module structure over the Cohomological Hall algebra (CoHa, for short). The major result of this section is Theorem \ref{thm_kernel} which states that the kernel of the quotient map from the CoHa to the Chow rings of non-commutative Hilbert schemes can be described using the CoHa-multiplication. The proof of this theorem relies on the Harder-Narasimhan stratification, when interpreting a non-commutative Hilbert scheme as a framed quiver moduli space. 

	\section*{Acknowledgements}

I would like to thank Markus Reineke and Sergey Mozgovoy for several very inspiring discussions concerning this work. I also wish to thank Yan Soibelman for sending me a preliminary version of his paper \cite{Soibelman:14} as well as Xinli Xiao for pointing out an imprecision in a previous version of the present article. The comments made by the referee have been great help in improving the exposition of this paper and were much appreciated.
At the time this paper was written, I was supported by the DFG priority program SPP 1388 ``Representation Theory''.
	\pagebreak[4]
	\section{Terminology and facts}
Fix an algebraically closed field $\kk$. We recollect the notation and some of the results of \cite{Reineke:05}.

\subsection{Notation}
\label{page_def_Rhat}

Fix positive integers $d, m$ and $n$ and vector spaces $V$ of dimension $n$ and $W$ of dimension $d$. Define $\smash{\hat{R}}$ to be the vector space $\Hom(V,W) \oplus \End(W)^m$ and let $G$ be the algebraic group $\Gl(W)$ which acts on $\smash{\hat{R}}$ via
$$
	g \cdot (f, \phi_1,\ldots,\phi_m) = (gf, g\phi_1 g^{-1}, \ldots, g\phi_m g^{-1}).
$$
An element $(f,\phi) = (f,\phi_1,\ldots,\phi_m)$ is called \textbf{stable} if $\kk\langle \phi_1,\ldots,\phi_m \rangle f(V) = W$, i.e. the image of $f$ generates $W$ regarded as a representation of the free non-commutative algebra $A = \kk\langle x_1,\ldots,x_m \rangle$ in $m$ variables. On the set $\smash{\hat{R}}^{\st}$ of stable points of $\smash{\hat{R}}$, a geometric $G$-quotient 
$$
	\pi: \hat{R}^{\st} \to \Hilb_{d,n}^{(m)}
$$
exists. It is even a principal $G$-bundle (by \cite[Prop.\ 0.9]{GIT:94}). The variety $\Hilb_{d,n}^{(m)}$ is called a \textbf{non-commutative Hilbert scheme}. As $m$ is fixed throughout this text, we sometimes suppress the dependency on $m$ and write $\Hilb_{d,n}$, for convenience. It is a smooth and irreducible variety of dimension $N := (m-1)d^2 + nd$. Its points parametrize $A$-submodules of codimension $d$ of the free module $A^n$. Denote the image $\pi(f,\phi)$ by $[f,\phi]$.

On $\hat{R}^{\st}$, we consider the $G$-equivariant vector bundle $\hat{R}^{\st} \times W \to \hat{R}^{\st}$, the trivial bundle, equipped with the $G$-action $g \cdot ((f,\phi),v) = ((gf, g\phi g^{-1}),gv)$. This descends to a vector bundle $\UU$ of rank $d$ on $\Hilb_{d,n}$, meaning there exists a vector bundle $\UU$ on $\Hilb_{d,n}$ such that $\pi^{-1}\UU$ with its canonical $G$-action equals the $G$-bundle we have just described. The $G$-linear endomorphisms $\smash{\hat{R}^{\st} \times W \to \hat{R}^{\st} \times W}$ mapping a point $((f,\phi),v)$ to $((f,\phi), \phi_iv)$ descend to endomorphisms $\Phi_i: \UU \to \UU$. Choosing a basis $e_1,\ldots,e_n$ of $V$ gives rise to $G$-linear sections $\smash{\hat{R}}^{\st} \to \smash{\hat{R}}^{\st} \times W$ defined by sending $(f,\phi)$ to $((f,\phi),fe_i)$. These sections, in turn, induce sections $s_1,\ldots,s_n$ of $\UU$.

The variety $\smash{\Hilb_{d,n}^{(m)}}$ is a so-called framed quiver moduli space (cf.\ \cite{ER:09}, \cite{Crawley-Boevey:03}, \cite{Nakajima:96}, or \cite{Reineke:08_Framed}). Consider the $m$-loop quiver $Q$ consisting of a single vertex $i$ and $m$ loops. A dimension vector for this quiver is just a natural number, say $d$. A representation of $Q$ of dimension vector $d$ is a $d$-dimensional vector space $W$ together with $m$ endomorphisms $\phi_1,\ldots,\phi_m$. In the sense of King \cite[Def.\ 2.1]{King:94} (which is a reformulation of Mumford's definition \cite[Def.\ 1.7]{GIT:94}), every representation of $Q$ is semi-stable, regardless of the choice of the linear form $\theta: \Z \to \Z$. Let's choose $\theta = 0$ for convenience. Define a quiver $\hat{Q}_n$ with two vertices $\infty$ and $i$ and with $n$ arrows pointing from $\infty$ to $i$ and $m$ loops at $i$. In a picture
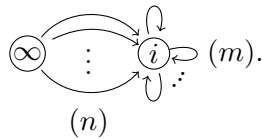
\begin{center}
	\begin{tikzpicture}[
	description/.style={fill=white,inner sep=2pt},
	implies/.style={double,double equal sign distance,-implies}
	]
	\matrix[column sep={3em}]
	{
		\node(a) {$\infty$}; & \node(b) {$i$}; \\
	};
	\path[->] (a) edge[bend left=50] (b)
		edge[bend left=30] (b)
		edge[bend right=50] (b);
	\node at ($(a) + (2.1em,0em)$) {$\vdots$};
	\path[->] (b) edge[loop above] ();
	\path[->] (b) edge[loop right] ();
	\path[->] (b) edge[loop below] ();
	\node at ($(b) + (1em,-.7em)$) {$.$};
	\node at ($(b) + (.9em,-.9em)$) {$.$};
	\node at ($(b) + (.7em,-1.1em)$) {$.$};
	\node at ($(b) + (2.8em,0em)$) {$(m).$};
	\node at ($(a) + (2.1em,-2.4em)$) {$(n)$};
	\draw (a) circle (.61em);
	\draw (b) circle (.54em);
\end{tikzpicture}
\end{center}
For this quiver, we consider the dimension vector $(1,d)$. A representation of $\hat{Q}_n$ consists of a representation of $Q$ of dimension vector $d$ and additionally, a linear map $f$ from an $n$-dimensional space $V$ to $W$. When choosing the extended stability condition $\smash{\hat{\theta}}$ as in \cite[Sect.\ 3]{Reineke:08_Framed} according to $\theta = 0$, we obtain that such a representation $(f,\phi)$ is stable if and only if it is $\smash{\hat{\theta}}$-semi-stable if and only if it is $\smash{\hat{\theta}}$-stable (see \cite[Prop.\ 3.3]{ER:09}).

\subsection{Words and forests}

Let $\Omega := \Omega^{(m)}$ be the set of words on the alphabet $\{1,\ldots,m\}$. The empty word will be denoted by $\epsilon$.

\begin{defn*}
	A finite subset $S$ of $\Omega$ is called a ($m$-\textbf{ary}) \textbf{tree} if it is closed under taking left subwords, that means, $w \in S$ provided $ww' \in S$ for some $w' \in \Omega$. A ($m$-\textbf{ary}) \textbf{forest} \textbf{with $n$ roots} is an $n$-tuple $S_* = (S_1,\ldots,S_n)$ of ($m$-ary) trees.
\end{defn*}

\begin{rem*}
	Note that the empty set is a tree. Thus, in a forest with, say, $n$ roots, we allow empty trees to occur.
\end{rem*}

Let $\FF_{d,n} := \FF_{d,n}^{(m)}$ be the set of $m$-ary forests with $n$ roots and $d$ nodes. Here, a forest $S_* = (S_1,\ldots,S_n)$ is said to have $d$ nodes if $\sharp S_1 + \cdots + \sharp S_n = d$. For a word $w$ with $w \in S_k$, we write $(k,w) \in S_*$.

A pair $(k',w')$ consisting of an index $1 \leq k' \leq n$ and a word $w' \in \Omega$ is called \textbf{critical} for a forest $S_*$ if either $w' = \epsilon$ and $S_{k'} = \emptyset$ or if $w' \notin S_{k'}$ but there exists a word $w \in S_{k'}$ and a letter $i \in \{1,\ldots,m\}$ with $w' = wi$. We define $C(S_*)$ to be the set of critical pairs of $S_*$. Its cardinality $c(S_*)$ equals $(m-1)\sharp S_* + n$.

We introduce an ordering on $\Omega$, the \textbf{lexicographic ordering}. For two words $w = i_1\cdots i_s$ and $w' = i'_1 \cdots i'_t$, let $p$ be the largest index such that $i_p = i'_p$. We formally define $p = 0$ if such an index doesn't exist. Define $w \leq w'$ if either $p = s$ (i.e.\ $w$ is a left subword of $w'$) or $i_{p+1} < i'_{p+1}$. This ordering can be extended to an ordering on the set of trees. Let $S$ and $S'$ be two distinct trees. Define $S < S'$ if either $\sharp S > \sharp S'$ or $\sharp S = \sharp S'$ and, writing $S = \{ w_1 < \cdots < w_s \}$ and $S' = \{ w'_1 < \cdots < w'_s \}$, we obtain $w_{p+1} < w'_{p+1}$ for the maximal index $p$ such that $w_p = w'_p$. Let's enlarge this to an ordering on the set $\FF_{d,n}$. For two distinct forests $S_*, S'_* \in \FF_{d,n}$, let $p$ be the largest index with $S_p = S'_p$ (again, $p = 0$ if $S_1 \neq S'_1$). Define $S_* < S'_*$ if $S_{p+1} < S'_{p+1}$.

Let $S_*$ be a forest. Define $D(S_*)$ to be the set of all quadruples $(k,w,k',w')$ consisting of indexes $1 \leq k,k' \leq n$ and words $w,w' \in \Omega$ satisfying
\begin{itemize}
	\item $(k,w) \in S_*$,
	\item $(k',w') \in C(S_*)$, and
	\item $(k,w) < (k',w')$, that means either $k < k'$ or $k = k'$ and $w < w'$.
\end{itemize}
The cardinality of $D(S_*)$ will be denoted $d(S_*)$.

\begin{ex*}
	Let's describe all $S_* \in \FF_{d,n}^{(m)}$ for $m = 2$, $d = 3$ and $n =1$. This example will accompany us throughout the text. When displaying $\Omega$ as follows
	
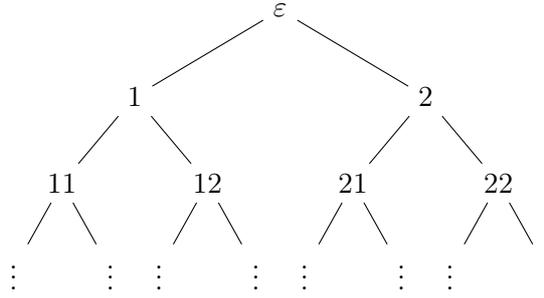
\begin{center}
		\begin{tikzpicture}[-,>=stealth',level/.style={sibling distance = 10em/#1, level distance = 3em}] 
			\node {$\epsilon$}
				child{ node {1} 
					child{ node {11}
						child { node {$\vdots$}}
						child { node {$\vdots$}}
					}
					child{ node {12}
						child { node {$\vdots$}}
						child { node {$\vdots$}}
					}                            
				}
				child{ node {2}
					child{ node {21}
						child { node {$\vdots$}}
						child { node {$\vdots$}}
					}
					child{ node {22}
						child { node {$\vdots$}}
						child { node {$\vdots$}}
					}
				}
			; 
		\end{tikzpicture}
	\end{center}
	then the $2$-ary - let's call them binary - trees with $3$ nodes are exactly those:
	
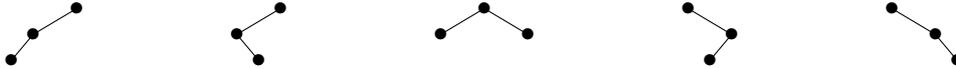
\begin{center}
		\begin{tikzpicture}[>=stealth',level/.style={sibling distance = 3em/#1, level distance = .9em}, parent anchor = center, child anchor = center] 
			\node (s) {}
				child { node (s1) {} edge from parent[black]
					child { node (s11) {} edge from parent[black]
						child { node (s111) {} edge from parent[white]}
						child { node (s112) {} edge from parent[white]}
					}
					child { node (s12) {} edge from parent[white]
						child { node (s121) {} edge from parent[white]}
						child { node (s122) {} edge from parent[white]}
					}
				}
				child { node (s2) {} edge from parent[white]
					child { node (s21) {} edge from parent[white]
						child { node (s211) {} edge from parent[white]}
						child { node (s212) {} edge from parent[white]}
					}
					child { node (s22) {} edge from parent[white]
						child { node (s221) {} edge from parent[white]}
						child { node (s222) {} edge from parent[white]}
					}
				}
			;
			\node [black] at (s) {$\bullet$};
			\node [black] at (s1) {$\bullet$};
			\node [black] at (s11) {$\bullet$};
			\node [white] at (s111) {$\bullet$};
			\node [white] at (s112) {$\bullet$};
			\node [white] at (s12) {$\bullet$};
			\node [white] at (s121) {$\bullet$};
			\node [white] at (s122) {$\bullet$};
			\node [white] at (s2) {$\bullet$};
			\node [white] at (s21) {$\bullet$};
			\node [white] at (s211) {$\bullet$};
			\node [white] at (s212) {$\bullet$};
			\node [white] at (s22) {$\bullet$};
			\node [white] at (s221) {$\bullet$};
			\node [white] at (s222) {$\bullet$};
		\end{tikzpicture}
		\begin{tikzpicture}[>=stealth',level/.style={sibling distance = 3em/#1, level distance = .9em}, parent anchor = center, child anchor = center] 
			\node (s) {}
				child { node (s1) {} edge from parent[black]
					child { node (s11) {} edge from parent[white]
						child { node (s111) {} edge from parent[white]}
						child { node (s112) {} edge from parent[white]}
					}
					child { node (s12) {} edge from parent[black]
						child { node (s121) {} edge from parent[white]}
						child { node (s122) {} edge from parent[white]}
					}
				}
				child { node (s2) {} edge from parent[white]
					child { node (s21) {} edge from parent[white]
						child { node (s211) {} edge from parent[white]}
						child { node (s212) {} edge from parent[white]}
					}
					child { node (s22) {} edge from parent[white]
						child { node (s221) {} edge from parent[white]}
						child { node (s222) {} edge from parent[white]}
					}
				}
			;
			\node [black] at (s) {$\bullet$};
			\node [black] at (s1) {$\bullet$};
			\node [white] at (s11) {$\bullet$};
			\node [white] at (s111) {$\bullet$};
			\node [white] at (s112) {$\bullet$};
			\node [black] at (s12) {$\bullet$};
			\node [white] at (s121) {$\bullet$};
			\node [white] at (s122) {$\bullet$};
			\node [white] at (s2) {$\bullet$};
			\node [white] at (s21) {$\bullet$};
			\node [white] at (s211) {$\bullet$};
			\node [white] at (s212) {$\bullet$};
			\node [white] at (s22) {$\bullet$};
			\node [white] at (s221) {$\bullet$};
			\node [white] at (s222) {$\bullet$};
		\end{tikzpicture}
		\begin{tikzpicture}[>=stealth',level/.style={sibling distance = 3em/#1, level distance = .9em}, parent anchor = center, child anchor = center] 
			\node (s) {}
				child { node (s1) {} edge from parent[black]
					child { node (s11) {} edge from parent[white]
						child { node (s111) {} edge from parent[white]}
						child { node (s112) {} edge from parent[white]}
					}
					child { node (s12) {} edge from parent[white]
						child { node (s121) {} edge from parent[white]}
						child { node (s122) {} edge from parent[white]}
					}
				}
				child { node (s2) {} edge from parent[black]
					child { node (s21) {} edge from parent[white]
						child { node (s211) {} edge from parent[white]}
						child { node (s212) {} edge from parent[white]}
					}
					child { node (s22) {} edge from parent[white]
						child { node (s221) {} edge from parent[white]}
						child { node (s222) {} edge from parent[white]}
					}
				}
			;
			\node [black] at (s) {$\bullet$};
			\node [black] at (s1) {$\bullet$};
			\node [white] at (s11) {$\bullet$};
			\node [white] at (s111) {$\bullet$};
			\node [white] at (s112) {$\bullet$};
			\node [white] at (s12) {$\bullet$};
			\node [white] at (s121) {$\bullet$};
			\node [white] at (s122) {$\bullet$};
			\node [black] at (s2) {$\bullet$};
			\node [white] at (s21) {$\bullet$};
			\node [white] at (s211) {$\bullet$};
			\node [white] at (s212) {$\bullet$};
			\node [white] at (s22) {$\bullet$};
			\node [white] at (s221) {$\bullet$};
			\node [white] at (s222) {$\bullet$};
		\end{tikzpicture}
		\begin{tikzpicture}[>=stealth',level/.style={sibling distance = 3em/#1, level distance = .9em}, parent anchor = center, child anchor = center] 
			\node (s) {}
				child { node (s1) {} edge from parent[white]
					child { node (s11) {} edge from parent[white]
						child { node (s111) {} edge from parent[white]}
						child { node (s112) {} edge from parent[white]}
					}
					child { node (s12) {} edge from parent[white]
						child { node (s121) {} edge from parent[white]}
						child { node (s122) {} edge from parent[white]}
					}
				}
				child { node (s2) {} edge from parent[black]
					child { node (s21) {} edge from parent[black]
						child { node (s211) {} edge from parent[white]}
						child { node (s212) {} edge from parent[white]}
					}
					child { node (s22) {} edge from parent[white]
						child { node (s221) {} edge from parent[white]}
						child { node (s222) {} edge from parent[white]}
					}
				}
			;
			\node [black] at (s) {$\bullet$};
			\node [white] at (s1) {$\bullet$};
			\node [white] at (s11) {$\bullet$};
			\node [white] at (s111) {$\bullet$};
			\node [white] at (s112) {$\bullet$};
			\node [white] at (s12) {$\bullet$};
			\node [white] at (s121) {$\bullet$};
			\node [white] at (s122) {$\bullet$};
			\node [black] at (s2) {$\bullet$};
			\node [black] at (s21) {$\bullet$};
			\node [white] at (s211) {$\bullet$};
			\node [white] at (s212) {$\bullet$};
			\node [white] at (s22) {$\bullet$};
			\node [white] at (s221) {$\bullet$};
			\node [white] at (s222) {$\bullet$};
		\end{tikzpicture}
		\begin{tikzpicture}[>=stealth',level/.style={sibling distance = 3em/#1, level distance = .9em}, parent anchor = center, child anchor = center] 
			\node (s) {}
				child { node (s1) {} edge from parent[white]
					child { node (s11) {} edge from parent[white]
						child { node (s111) {} edge from parent[white]}
						child { node (s112) {} edge from parent[white]}
					}
					child { node (s12) {} edge from parent[white]
						child { node (s121) {} edge from parent[white]}
						child { node (s122) {} edge from parent[white]}
					}
				}
				child { node (s2) {} edge from parent[black]
					child { node (s21) {} edge from parent[white]
						child { node (s211) {} edge from parent[white]}
						child { node (s212) {} edge from parent[white]}
					}
					child { node (s22) {} edge from parent[black]
						child { node (s221) {} edge from parent[white]}
						child { node (s222) {} edge from parent[white]}
					}
				}
			;
			\node [black] at (s) {$\bullet$};
			\node [white] at (s1) {$\bullet$};
			\node [white] at (s11) {$\bullet$};
			\node [white] at (s111) {$\bullet$};
			\node [white] at (s112) {$\bullet$};
			\node [white] at (s12) {$\bullet$};
			\node [white] at (s121) {$\bullet$};
			\node [white] at (s122) {$\bullet$};
			\node [black] at (s2) {$\bullet$};
			\node [white] at (s21) {$\bullet$};
			\node [white] at (s211) {$\bullet$};
			\node [white] at (s212) {$\bullet$};
			\node [black] at (s22) {$\bullet$};
			\node [white] at (s221) {$\bullet$};
			\node [white] at (s222) {$\bullet$};
		\end{tikzpicture}
	\end{center}
	Compare this to Stanley's list of descriptions of the Catalan numbers (cf.\ \cite[Ex. 6.19]{Stanley:99}). The above is part (c) of the list. When considering the sets $S_* \sqcup C(S_*)$, we also get a tree (or a forest, in general), more precisely a \emph{plane} binary tree with $(m-1)d+n = 7$ vertices. This is part (d) of Stanley's list. In the sketch below, the vertices belonging to $S_*$ are displayed in gray:
	
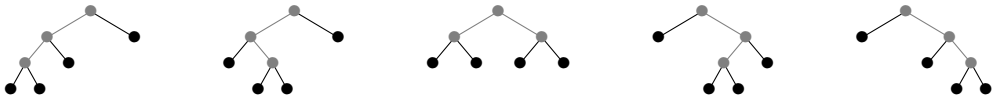
\begin{center}
		\begin{tikzpicture}[>=stealth',level/.style={sibling distance = 3em/#1, level distance = .9em}, parent anchor = center, child anchor = center] 
			\node (s) {}
				child { node (s1) {} edge from parent[gray]
					child { node (s11) {} edge from parent[gray]
						child { node (s111) {} edge from parent[black]}
						child { node (s112) {} edge from parent[black]}
					}
					child { node (s12) {} edge from parent[black]
						child { node (s121) {} edge from parent[white]}
						child { node (s122) {} edge from parent[white]}
					}
				}
				child { node (s2) {} edge from parent[black]
					child { node (s21) {} edge from parent[white]
						child { node (s211) {} edge from parent[white]}
						child { node (s212) {} edge from parent[white]}
					}
					child { node (s22) {} edge from parent[white]
						child { node (s221) {} edge from parent[white]}
						child { node (s222) {} edge from parent[white]}
					}
				}
			;
			\node [gray] at (s) {$\bullet$};
			\node [gray] at (s1) {$\bullet$};
			\node [gray] at (s11) {$\bullet$};
			\node [black] at (s111) {$\bullet$};
			\node [black] at (s112) {$\bullet$};
			\node [black] at (s12) {$\bullet$};
			\node [white] at (s121) {$\bullet$};
			\node [white] at (s122) {$\bullet$};
			\node [black] at (s2) {$\bullet$};
			\node [white] at (s21) {$\bullet$};
			\node [white] at (s211) {$\bullet$};
			\node [white] at (s212) {$\bullet$};
			\node [white] at (s22) {$\bullet$};
			\node [white] at (s221) {$\bullet$};
			\node [white] at (s222) {$\bullet$};
		\end{tikzpicture}
		\begin{tikzpicture}[>=stealth',level/.style={sibling distance = 3em/#1, level distance = .9em}, parent anchor = center, child anchor = center] 
			\node (s) {}
				child { node (s1) {} edge from parent[gray]
					child { node (s11) {} edge from parent[black]
						child { node (s111) {} edge from parent[white]}
						child { node (s112) {} edge from parent[white]}
					}
					child { node (s12) {} edge from parent[gray]
						child { node (s121) {} edge from parent[black]}
						child { node (s122) {} edge from parent[black]}
					}
				}
				child { node (s2) {} edge from parent[black]
					child { node (s21) {} edge from parent[white]
						child { node (s211) {} edge from parent[white]}
						child { node (s212) {} edge from parent[white]}
					}
					child { node (s22) {} edge from parent[white]
						child { node (s221) {} edge from parent[white]}
						child { node (s222) {} edge from parent[white]}
					}
				}
			;
			\node [gray] at (s) {$\bullet$};
			\node [gray] at (s1) {$\bullet$};
			\node [black] at (s11) {$\bullet$};
			\node [white] at (s111) {$\bullet$};
			\node [white] at (s112) {$\bullet$};
			\node [gray] at (s12) {$\bullet$};
			\node [black] at (s121) {$\bullet$};
			\node [black] at (s122) {$\bullet$};
			\node [black] at (s2) {$\bullet$};
			\node [white] at (s21) {$\bullet$};
			\node [white] at (s211) {$\bullet$};
			\node [white] at (s212) {$\bullet$};
			\node [white] at (s22) {$\bullet$};
			\node [white] at (s221) {$\bullet$};
			\node [white] at (s222) {$\bullet$};
		\end{tikzpicture}
		\begin{tikzpicture}[>=stealth',level/.style={sibling distance = 3em/#1, level distance = .9em}, parent anchor = center, child anchor = center] 
			\node (s) {}
				child { node (s1) {} edge from parent[gray]
					child { node (s11) {} edge from parent[black]
						child { node (s111) {} edge from parent[white]}
						child { node (s112) {} edge from parent[white]}
					}
					child { node (s12) {} edge from parent[black]
						child { node (s121) {} edge from parent[white]}
						child { node (s122) {} edge from parent[white]}
					}
				}
				child { node (s2) {} edge from parent[gray]
					child { node (s21) {} edge from parent[black]
						child { node (s211) {} edge from parent[white]}
						child { node (s212) {} edge from parent[white]}
					}
					child { node (s22) {} edge from parent[black]
						child { node (s221) {} edge from parent[white]}
						child { node (s222) {} edge from parent[white]}
					}
				}
			;
			\node [gray] at (s) {$\bullet$};
			\node [gray] at (s1) {$\bullet$};
			\node [black] at (s11) {$\bullet$};
			\node [white] at (s111) {$\bullet$};
			\node [white] at (s112) {$\bullet$};
			\node [black] at (s12) {$\bullet$};
			\node [white] at (s121) {$\bullet$};
			\node [white] at (s122) {$\bullet$};
			\node [gray] at (s2) {$\bullet$};
			\node [black] at (s21) {$\bullet$};
			\node [white] at (s211) {$\bullet$};
			\node [white] at (s212) {$\bullet$};
			\node [black] at (s22) {$\bullet$};
			\node [white] at (s221) {$\bullet$};
			\node [white] at (s222) {$\bullet$};
		\end{tikzpicture}
		\begin{tikzpicture}[>=stealth',level/.style={sibling distance = 3em/#1, level distance = .9em}, parent anchor = center, child anchor = center] 
			\node (s) {}
				child { node (s1) {} edge from parent[black]
					child { node (s11) {} edge from parent[white]
						child { node (s111) {} edge from parent[white]}
						child { node (s112) {} edge from parent[white]}
					}
					child { node (s12) {} edge from parent[white]
						child { node (s121) {} edge from parent[white]}
						child { node (s122) {} edge from parent[white]}
					}
				}
				child { node (s2) {} edge from parent[gray]
					child { node (s21) {} edge from parent[gray]
						child { node (s211) {} edge from parent[black]}
						child { node (s212) {} edge from parent[black]}
					}
					child { node (s22) {} edge from parent[black]
						child { node (s221) {} edge from parent[white]}
						child { node (s222) {} edge from parent[white]}
					}
				}
			;
			\node [gray] at (s) {$\bullet$};
			\node [black] at (s1) {$\bullet$};
			\node [white] at (s11) {$\bullet$};
			\node [white] at (s111) {$\bullet$};
			\node [white] at (s112) {$\bullet$};
			\node [white] at (s12) {$\bullet$};
			\node [white] at (s121) {$\bullet$};
			\node [white] at (s122) {$\bullet$};
			\node [gray] at (s2) {$\bullet$};
			\node [gray] at (s21) {$\bullet$};
			\node [black] at (s211) {$\bullet$};
			\node [black] at (s212) {$\bullet$};
			\node [black] at (s22) {$\bullet$};
			\node [white] at (s221) {$\bullet$};
			\node [white] at (s222) {$\bullet$};
		\end{tikzpicture}
		\begin{tikzpicture}[>=stealth',level/.style={sibling distance = 3em/#1, level distance = .9em}, parent anchor = center, child anchor = center] 
			\node (s) {}
				child { node (s1) {} edge from parent[black]
					child { node (s11) {} edge from parent[white]
						child { node (s111) {} edge from parent[white]}
						child { node (s112) {} edge from parent[white]}
					}
					child { node (s12) {} edge from parent[white]
						child { node (s121) {} edge from parent[white]}
						child { node (s122) {} edge from parent[white]}
					}
				}
				child { node (s2) {} edge from parent[gray]
					child { node (s21) {} edge from parent[black]
						child { node (s211) {} edge from parent[white]}
						child { node (s212) {} edge from parent[white]}
					}
					child { node (s22) {} edge from parent[gray]
						child { node (s221) {} edge from parent[black]}
						child { node (s222) {} edge from parent[black]}
					}
				}
			;
			\node [gray] at (s) {$\bullet$};
			\node [black] at (s1) {$\bullet$};
			\node [white] at (s11) {$\bullet$};
			\node [white] at (s111) {$\bullet$};
			\node [white] at (s112) {$\bullet$};
			\node [white] at (s12) {$\bullet$};
			\node [white] at (s121) {$\bullet$};
			\node [white] at (s122) {$\bullet$};
			\node [gray] at (s2) {$\bullet$};
			\node [black] at (s21) {$\bullet$};
			\node [white] at (s211) {$\bullet$};
			\node [white] at (s212) {$\bullet$};
			\node [gray] at (s22) {$\bullet$};
			\node [black] at (s221) {$\bullet$};
			\node [black] at (s222) {$\bullet$};
		\end{tikzpicture}
	\end{center}
	Now, let's determine the sets $D(S_*)$. They are given by
	
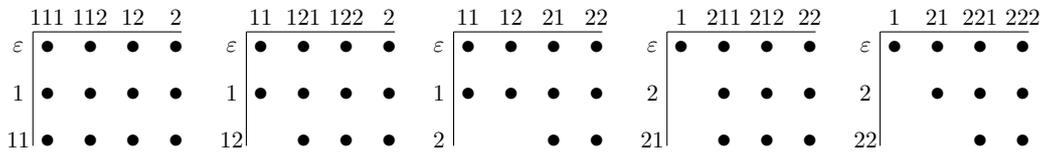
\begin{center}
		\begin{tikzpicture}[description/.style={fill=white,inner sep=2pt}]
			\matrix(m)[matrix of math nodes, row sep=.2em, column sep=.3em, text height=1.5ex, text depth=0.25ex]
			{
				\bullet		& \bullet	& \bullet	& \bullet \\
				\bullet		& \bullet	& \bullet	& \bullet \\
				\bullet		& \bullet	& \bullet	& \bullet \\
			};
			\node[scale=.8] at ($(m-1-1)+(-1em,0em)$) {$\epsilon$};
			\node[scale=.8] at ($(m-2-1)+(-1em,0em)$) {$1$};
			\node[scale=.8] at ($(m-3-1)+(-1em,0em)$) {$11$};
			\node[scale=.8] at ($(m-1-1)+(0em,1em)$) {$111$};
			\node[scale=.8] at ($(m-1-2)+(0em,1em)$) {$112$};
			\node[scale=.8] at ($(m-1-3)+(0em,1em)$) {$12$};
			\node[scale=.8] at ($(m-1-4)+(0em,1em)$) {$2$};
			\draw ($(m-1-1)+(-.5em,.5em)$) -- ($(m-3-1)+(-.5em,-.2em)$);
			\draw ($(m-1-1)+(-.5em,.5em)$) -- ($(m-1-4)+(.2em,.5em)$);
		\end{tikzpicture}
		\begin{tikzpicture}[description/.style={fill=white,inner sep=2pt}]
			\matrix(m)[matrix of math nodes, row sep=.2em, column sep=.3em, text height=1.5ex, text depth=0.25ex]
			{
				\bullet		& \bullet	& \bullet	& \bullet \\
				\bullet		& \bullet	& \bullet	& \bullet \\
						& \bullet	& \bullet	& \bullet \\
			};
			\node[scale=.8] at ($(m-1-1)+(-1em,0em)$) {$\epsilon$};
			\node[scale=.8] at ($(m-2-1)+(-1em,0em)$) {$1$};
			\node[scale=.8] at ($(m-3-1)+(-1em,0em)$) {$12$};
			\node[scale=.8] at ($(m-1-1)+(0em,1em)$) {$11$};
			\node[scale=.8] at ($(m-1-2)+(0em,1em)$) {$121$};
			\node[scale=.8] at ($(m-1-3)+(0em,1em)$) {$122$};
			\node[scale=.8] at ($(m-1-4)+(0em,1em)$) {$2$};
			\draw ($(m-1-1)+(-.5em,.5em)$) -- ($(m-3-1)+(-.5em,-.2em)$);
			\draw ($(m-1-1)+(-.5em,.5em)$) -- ($(m-1-4)+(.2em,.5em)$);
		\end{tikzpicture}
		\begin{tikzpicture}[description/.style={fill=white,inner sep=2pt}]
			\matrix(m)[matrix of math nodes, row sep=.2em, column sep=.3em, text height=1.5ex, text depth=0.25ex]
			{
				\bullet		& \bullet	& \bullet	& \bullet \\
				\bullet		& \bullet	& \bullet	& \bullet \\
						& 		& \bullet	& \bullet \\
			};
			\node[scale=.8] at ($(m-1-1)+(-1em,0em)$) {$\epsilon$};
			\node[scale=.8] at ($(m-2-1)+(-1em,0em)$) {$1$};
			\node[scale=.8] at ($(m-3-1)+(-1em,0em)$) {$2$};
			\node[scale=.8] at ($(m-1-1)+(0em,1em)$) {$11$};
			\node[scale=.8] at ($(m-1-2)+(0em,1em)$) {$12$};
			\node[scale=.8] at ($(m-1-3)+(0em,1em)$) {$21$};
			\node[scale=.8] at ($(m-1-4)+(0em,1em)$) {$22$};
			\draw ($(m-1-1)+(-.5em,.5em)$) -- ($(m-3-1)+(-.5em,-.2em)$);
			\draw ($(m-1-1)+(-.5em,.5em)$) -- ($(m-1-4)+(.2em,.5em)$);
		\end{tikzpicture}
		\begin{tikzpicture}[description/.style={fill=white,inner sep=2pt}]
			\matrix(m)[matrix of math nodes, row sep=.2em, column sep=.3em, text height=1.5ex, text depth=0.25ex]
			{
				\bullet		& \bullet	& \bullet	& \bullet \\
						& \bullet	& \bullet	& \bullet \\
						& \bullet	& \bullet	& \bullet \\
			};
			\node[scale=.8] at ($(m-1-1)+(-1em,0em)$) {$\epsilon$};
			\node[scale=.8] at ($(m-2-1)+(-1em,0em)$) {$2$};
			\node[scale=.8] at ($(m-3-1)+(-1em,0em)$) {$21$};
			\node[scale=.8] at ($(m-1-1)+(0em,1em)$) {$1$};
			\node[scale=.8] at ($(m-1-2)+(0em,1em)$) {$211$};
			\node[scale=.8] at ($(m-1-3)+(0em,1em)$) {$212$};
			\node[scale=.8] at ($(m-1-4)+(0em,1em)$) {$22$};
			\draw ($(m-1-1)+(-.5em,.5em)$) -- ($(m-3-1)+(-.5em,-.2em)$);
			\draw ($(m-1-1)+(-.5em,.5em)$) -- ($(m-1-4)+(.2em,.5em)$);
		\end{tikzpicture}
		\begin{tikzpicture}[description/.style={fill=white,inner sep=2pt}]
			\matrix(m)[matrix of math nodes, row sep=.2em, column sep=.3em, text height=1.5ex, text depth=0.25ex]
			{
				\bullet		& \bullet	& \bullet	& \bullet \\
						& \bullet	& \bullet	& \bullet \\
						& 		& \bullet	& \bullet \\
			};
			\node[scale=.8] at ($(m-1-1)+(-1em,0em)$) {$\epsilon$};
			\node[scale=.8] at ($(m-2-1)+(-1em,0em)$) {$2$};
			\node[scale=.8] at ($(m-3-1)+(-1em,0em)$) {$22$};
			\node[scale=.8] at ($(m-1-1)+(0em,1em)$) {$1$};
			\node[scale=.8] at ($(m-1-2)+(0em,1em)$) {$21$};
			\node[scale=.8] at ($(m-1-3)+(0em,1em)$) {$221$};
			\node[scale=.8] at ($(m-1-4)+(0em,1em)$) {$222$};
			\draw ($(m-1-1)+(-.5em,.5em)$) -- ($(m-3-1)+(-.5em,-.2em)$);
			\draw ($(m-1-1)+(-.5em,.5em)$) -- ($(m-1-4)+(.2em,.5em)$);
		\end{tikzpicture}
	\end{center}
	and when viewing the missing entries as Young diagrams (after turning them upside down) fitting in some triangular shape, we obtain (vv) of \cite[Ex.~6.19]{Stanley:99}.
\end{ex*}

\subsection{A cell decomposition}

For a word $w \in \Omega$, say $w = i_1\cdots i_s$, and a point $(f,\phi) \in \hat{R}$, define the endomorphism $\phi_w$ of $W$ to be the composition $\phi_{i_s}\cdots\phi_{i_1}$. In the same vein, define $\Phi_w := \Phi_{i_s}\cdots \Phi_{i_1}$ to obtain an endomorphism of the bundle $\UU$. Finally, define the section $s_{(k,w)}$ of $\UU$ to be $\Phi_w s_k$.

\begin{defn*}
	Let $S_* \in \FF_{d,n}^{(m)}$ be a forest. Define $U_{S_*}$ to be the subset of all $[f,\phi] \in \smash{\Hilb_{d,n}^{(m)}}$ such that the vectors $\phi_wfe_k$ with $(k,w) \in S_*$ form a basis of~$W$.
\end{defn*}

Reineke shows in \cite[Le.\ 3.2]{Reineke:05} that for every point $[f,\phi]$ of $\Hilb_{d,n}$ and every forest $S'_*$ for which the tuple of vectors $(\phi_wfe_k \mid (k,w) \in S'_*)$ is linearly independent, there exists a forest $S_* \in \FF_{d,n}$ containing $S'_*$ such that $[f,\phi]$ is in $U_{S_*}$. Furthermore, by expressing $\phi_{w'}fe_{k'}$ in terms of the basis $\phi_wfe_k$ with $(k,w,k',w') \in D(S_*)$, he shows that $U_{S_*}$ is isomorphic to an affine space. This implies (cf.\ \cite[Cor.\ 3.3, Le.\ 3.4]{Reineke:05}) that the variety $\Hilb_{d,n}$ is covered by the open subsets $U_{S_*}$ with $S_* \in \FF_{d,n}$, each of which is isomorphic to an affine space of dimension $N = (m-1)d^2 + nd$.

Next, we define certain closed subsets of the $U_{S_*}$. These subsets will be the cells of the cell decomposition we are about to describe.

\begin{defn*}
	Let $S_* \in \smash{\FF_{d,n}^{(m)}}$ be a forest. Define $Z_{S_*}$ to be the set of all $[f,\phi] \in U_{S_*}$ such that for all critical pairs $(k',w') \in C(S_*)$, the vector $\phi_{w'}fe_{k'}$ is contained in the span of all $\phi_wfe_k$ with $(k,w) \in S$ and $(k,w) < (k',w')$.
\end{defn*}

In \cite{Reineke:05}, a description of $Z_{S_*}$ as a set in terms of the $U_{S'_*}$ for $S'_* < S_*$ is given. It reads as follows:

\begin{thm}[{\cite[Thm.\ 3.6]{Reineke:05}}]
	For all forests $S_* \in \FF_{d,n}^{(m)}$, we obtain $$Z_{S_*} = U_{S_*} \setminus \bigcup_{S'_* < S_*} U_{S'_*}.$$
\end{thm}

%
Moreover, equipping $Z_{S_*}$ with the reduced closed subscheme structure of $U_{S_*}$ and displaying $[f,\phi]$ in terms of the basis $\phi_wfe_k$ with $(k,w) \in S_*$, \cite[Le.~3.9]{Reineke:05} shows that $Z_{S_*}$ is isomorphic to an affine space of dimension $d(S_*)$.
%
%
These results lead to a main result of \cite{Reineke:05}, the existence of a cell decomposition of $\Hilb_{d,n}$. By definition, a cell decomposition of a variety is a descending sequence of closed subsets such that the successive complements are isomorphic to affine spaces. Define $A_{S_*} := \Hilb_{d,n} - \bigcup_{S'_* < S_*} U_{S'_*}$. Enumerating the forests of $\FF_{d,n}$ lexicographically, say $S_*^1 < \cdots < S_*^u$, and abbreviating $A_i := A_{S_*^i}$, we obtain a filtration
$$
	\Hilb_{d,n} = A_1 \supseteq A_2 \supseteq \cdots \supseteq A_u \supseteq A_{u+1} := \emptyset
$$
by closed subsets satisfying $A_i - A_{i+1} = Z_{S_*^i}$. Cutting a long story short:

\begin{thm}[{\cite[Thm.\ 1.3]{Reineke:05}}] \label{thm_cell}
	The non-commutative Hilbert scheme $\Hilb_{d,n}^{(m)}$ possesses a cell decomposition parametrized by forests $S_* \in \smash{\FF_{d,n}^{(m)}}$, whose cells $Z_{S_*}$ are of dimensions $d(S_*)$.
\end{thm}

An immediate application yields a basis of the Chow group of $\Hilb_{d,n}$. We denote by $\ZZ_{S_*}$ the closure of $Z_{S_*}$ in $\Hilb_{d,n}$ equipped with the reduced closed subscheme structure. As $\ZZ_{S_*}$ is irreducible, it becomes a closed subvariety of $\Hilb_{d,n}$.

\begin{cor}[{\cite[Cor.\ 4.3]{Reineke:05}}] \label{free}
	The Chow group $A_*(\Hilb_{d,n}^{(m)})$ is the free abelian group with basis $[\ZZ_{S_*}]$ for $S_* \in \smash{\FF_{d,n}^{(m)}}$.
\end{cor}

\begin{ex*}[continued]
	Again, let $m=2$, $d=3$ and $n=1$. We describe the cells $Z_{S_*}$ belonging to the binary trees $S_*$ with $3$ nodes. A point of $R$ may be viewed as a triple $(v,A,B)$, where $v \in \kk^3$ and $A$ and $B$ are $(3 \times 3)$-matrices. Write $[v,A,B]$ for its image in the non-commutative Hilbert scheme. The cells are
	\begin{align*}
		Z_{\treeA} &= \left\{ [v,A,B] \mid v, Av, A^2v \text{ basis of } \kk^3 \right\}, \\
		Z_{\treeB} &= \left\{ [v,A,B] \mid v, Av, BAv \text{ basis of } \kk^3 \text{ and } A^2v \in \langle v,Av \rangle \right\}, \\
		Z_{\treeC} &= \left\{ [v,A,B] \mid v, Av, Bv \text{ basis of } \kk^3 \text{ and } A^2v, BAv \in \langle v,Av \rangle \right\}, \\
		Z_{\treeD} &= \left\{ [v,A,B] \mid v, Bv, ABv \text{ basis of } \kk^3 \text{ and } Av \in \langle v \rangle \right\}, \text{ and} \\
		Z_{\treeE} &= \left\{ [v,A,B] \mid v, Bv, B^2v \text{ basis of } \kk^3,\ Av \in \langle v \rangle \text{, and } ABv \in \langle v,Bv \rangle \right\}.
	\end{align*}
	Knowing their dimensions, we are able to determine the Poincar\'{e} polynomial 
	$$
		\sum_i \dim_\Q(A_i(\Hilb_{d,n}^{(m)})_\Q) t^i
	$$ 
	at once. It reads $t^{12} + t^{11} + 2t^{10} + t^{9}$.
\end{ex*}

	\section{Another basis of the Chow group}

We are interested in the ring structure on the Chow group of the non-commutative Hilbert scheme $\smash{\Hilb_{d,n}^{(m)}}$. It turns out that computing the intersection product of two cell closures is rather difficult. We would therefore like to find another basis that is better adapted to the multiplication. This basis will be provided by Chern classes of the universal bundle $\UU$ which we have already introduced in the previous section.

\subsection{A connection between cell closures and Chern classes of the universal bundle}

Let $S_* \in \FF_{d,n}$ be a forest. Order the words of the trees lexicographically, i.e. $S_k = \{w_{k,1} < \cdots < w_{k,d_k} \}$. Consider all pairs $(k,w) \in S_*$ and order them lexicographically, too. This gives
$$
	(1,w_{1,1}) < \cdots < (1,w_{1,d_1}) < \cdots  < (n, w_{n,1}) < \cdots < (n,w_{n,d_n})
$$
and we denote these pairs as $x_1 < \cdots < x_d$. This means $(k,w_{k,\nu}) = x_{d_1 + \cdots + d_{k-1} + \nu}$. For a critical pair $x' = (k',w')$ of $S_*$, let $j = j_{S_*}(x')$ be the maximal index such that $x_j < x'$. Formally, let $j = 0$ if such an index doesn't exist. We express $j$ in a slightly different way. If $w' = \epsilon$, then $j = d_1 + \cdots + d_{k'-1}$. Otherwise, it is $j = d_1 + \cdots + d_{k'-1} + \nu$, where $\nu$ is the maximal index such that $w_\nu < w'$ (and $\nu$ is not zero in this case, but possibly $d_{k'}$). As $D(S_*)$ is clearly in bijection to the disjoint union of the sets $\{ (x_1,x'),\ldots,(x_{j(x')},x') \}$, with $x'$ ranging over all critical pairs of $S_*$, we see that $\sum_{x' \in C(S_*)} j(x') = d(S_*)$. Define $i(x') := i_{S_*}(x')$ to be $d - j(x')$. We will show the following:

\begin{thm} \label{thm}
	For all forests $S_* \in \FF_{d,n}^{(m)}$, we have
	$$
		\prod_{x' \in C(S_*)} c_{i(x')}(\UU) \cap [\Hilb_{d,n}^{(m)}]\ =\ [\ZZ_{S_*}]\ + \sum_{\substack{S''_* > S_*,\\ d(S''_*) = d(S_*)}}\!\!\!\! n_{S_*,S''_*} [\ZZ_{S''_*}]
	$$
	for some positive integers $n_{S_*,S''_*}$.
\end{thm}

Recall the section $s_x = s_{(k,w)}$ associated to any pair $x = (k,w)$ consisting of an index $1 \leq k \leq n$ and a word $w$. For a forest $S_*$, define $D_{S_*}$ as the intersection of the degeneracy loci
$$
	D_{S_*} = \bigcap_{x' \in C(S_*)} D_{S_*}(x'),
$$
where $D_{S_*}(x') := D(s_{x_1},\ldots,s_{x_{j(x')}},s_{x'})$ is the degeneracy locus as defined in \cite[Chap.\ 14]{Fulton:98}. As all degeneracy loci possess a natural structure of a closed subscheme of $\Hilb_{d,n}$, the subset $D_{S_*}$ does, too.

\begin{lem} \label{underl_closed_subset}
	The underlying closed subset of $D_{S_*}$ is 
    $$A_{S_*} = \Hilb_{d,n}^{(m)} - \bigcup_{S'_* < S_*} U_{S'_*}.$$
\end{lem}

\begin{proof}
	Order the pairs of $S_*$, i.e. $x_1 = (k_1,w_1) < \cdots < x_d = (k_d,w_d)$. A point $[f,\phi]$ lies in $D_{S_*}$ if and only if the vectors
	$$
		\phi_{w_1}fe_{k_1},\ldots,\phi_{w_{j(x')}}fe_{k_{j(x')}},\phi_{w'}fe_{k'}
	$$
	are linearly dependent for all critical pairs $x' = (k',w')$ of $S_*$. Let $S'_* \in \FF$ with $S'_* < S_*$, say $S' = \{x'_1 < \cdots < x'_d\}$. Define $p$ to be the maximal index such that $x_p = x'_p$. We have $x'_{p+1} < x_{p+1}$ and therefore $x'_{p+1} \notin S_*$. We write $x'_{p+1}$ as $(k',w')$, this means $w' \notin S_{k'}$. If $w' = \epsilon$, then $S_{k'} = \emptyset$ and if $w'$ is not the empty word, we write $w' = wi$ for some $w \in S'_{k'}$. As $(k',w) < x'_{p+1}$, we obtain $w \in S_{k'}$. This means $x'_{p+1}$ is a critical pair for $S_*$. Moreover, we get $j(x'_{p+1}) = p$. Therefore, the vectors
	$$
		\phi_{w_1}fe_{k_1},\ldots,\phi_{w_p}fe_p,\phi_{w'}fe_{k'}
	$$
	are linearly dependent and this implies that $[f,\phi]$ is not contained in $U_{S'_*}$. Conversely, assume that $[f,\phi]$ does not belong to the union $\bigcup_{S'_* < S_*} U_{S'_*}$. Let $x' = (k',w')$ be a critical pair for $S_*$. Let $j := j(x')$ and write $x_j = (k',w_{k',\nu})$. Consider the forest $S'_*$ consisting of $S'_k := S_k$ for all $k < k'$, of
	$$
		S'_{k'} = \{ w_{k',1} < \cdots < w_{k',\nu} < w' \}
	$$
	and of $S'_k = \emptyset$ for $k > k'$. Assume that the vectors $\phi_{w_1}fe_{k_1},\ldots,\phi_{w_j}fe_{k_j},\phi_{w'}fe_{k'}$ were linearly independent. By \cite[Le.\ 3.2]{Reineke:05}, there exists a forest $S''_*$ containing $S'_*$ such that $[f,\phi]$ belongs to $U_{S''_*}$. But this forest fulfills $S''_* < S_*$. A contradiction.
\end{proof}

\begin{ex*}[continued]
	Let's determine the underlying closed subsets $A_{S_*}$ of the $D_{S_*}$ in the - by now well known - case $m=2$, $d=3$ and $n=1$. We have
	\begin{align*}
		A_{\treeA} &= \Hilb_{3,1}^{(2)}, \\
		A_{\treeB} &= \left\{ [v,A,B] \mid v,Av,A^2v \text{ linearly dependent} \right\}, \\
		A_{\treeC} &= \left\{ [v,A,B] \mid v,Av,A^2v \text{ and } v,Av,BAv \text{ both linearly dependent} \right\}, \\
		A_{\treeD} &= \left\{ [v,A,B] \mid v,Av \text{ linearly dependent} \right\}, \text{ and}\\
		A_{\treeE} &= \left\{ [v,A,B] \mid v,Av \text{ and } v,Bv,ABv \text{ both linearly dependent} \right\}.
	\end{align*}
	We can easily see that the successive complements are, indeed, the cells $Z_{S_*}$.
\end{ex*}

For general reasons (cf.\ \cite[Thm.\ 14.4]{Fulton:98}), we know that every irreducible component of $D_{S_*}(x')$ has dimension at least $N - i(x')$. We will show that, in fact, equality holds.

\begin{lem} \label{lem_codim}
	Let $T_* = \{x_1 < \cdots < x_j < x' \} \in \FF_{j+1,n}^{(m)}$ be a forest. Then, the closed subset $D_{T_*} = D(s_{x_1},\ldots,s_{x_j},s_{x'})$ has pure dimension $N - (d-j)$ (or is empty).
\end{lem}

\begin{proof}
	The proof proceeds by induction on $j$. In the case $j = 0$, the forest $T_*$ equals $\{ (k',\epsilon) \}$ for an index $1 \leq k' \leq n$. Choose a forest $S_* \in \FF_{d,n}$ such that $(k',\epsilon) \in C(S_*)$. If $n=1$, such a forest does not exist and $D_{T_*}$ is empty. Otherwise, $D_{T_*} \cap U_{S_*} \neq \emptyset$. We apply the isomorphism $U_{S_*} \cong \A^N$ from \cite{Reineke:05}, provided by the functions $\lambda_{x,x'}$ on $U_{S_*}$ for every $x' = (k',w') \in C(S_*)$ and every $x \in S_*$. By definition, $\lambda_{x,x'}[f,\phi]$ is the coefficient occurring in the linear combination
	$$
		\phi_{w'}fe_{k'} = \sum_{x = (k,w) \in S_*} \lambda_{x,x'}[f,\phi] \cdot \phi_wfe_k
	$$
	for $x' = (k',w') \in C(S_*)$. The closed subscheme $D_{T_*} \cap U_{S_*}$ is defined by annihilation of all functions $\lambda_{x,(k',\epsilon)}$. This describes an affine space of dimension $N - d$.
	Assume that $j > 0$. Let $T'_* := \{x_1 < \cdots < x_j \}$. This is also a forest. Let $S_* \in \FF_{d,n}$ be a forest which contains $T'_*$. We consider $D_{T_*} \cap U_{S_*}$. If $x' \in S_*$ then $D_{T_*} \cap U_{S_*}$ is empty. So we require $x' \notin S_*$ which implies that $x'$ is a critical pair for $S_*$. Via the isomorphism $U_{S_*} \cong \A^N$ described above, $D_{T_*} \cap U_{S_*}$ is given by the ideal generated by all functions $\lambda_{x,x'}$ with $x \in S_*$ and $x > x'$. This describes an affine space of dimension $N - (d-j)$. This shows that every irreducible component of $D_{T_*}$ which intersects the open subset $U$, defined as the union $\bigcup_{S_*\supseteq T'_*} U_{S_*}$, has dimension $N-(d-j)$. To conclude the proof, we show that there are no irreducible components of $D_{T_*}$ which are contained in the complement of $U$. Assume there were such a component. As the intersection $D_{T_*} \cap U^c$ lies in $D_{T'_*}$, this component would have to be contained in an irreducible component of $D_{T'_*}$. But by induction hypothesis, we know that $D_{T'_*}$ has pure dimension $N-(d-j)-1$ which contradicts the fact that every irreducible component of $D_{T_*}$ has dimension at least $N-(d-j)$.
\end{proof}


The above lemma implies, using \cite[Ex.\ 14.4.2]{Fulton:98}, that the cycle $[D_{S_*}(x')]$ associated to the degeneracy locus $D_{S_*}(x') = D(s_{x_1},\ldots,s_{x_{j(x')}},s_{x'})$, regarded as an element of $A_{N-i(x')}(\Hilb_{d,n})$, equals $c_{i(x')}(\UU) \cap [\Hilb_{d,n}]$. We use this observation to prove Theorem \ref{thm}.

\begin{proof}[\normalfont \textit{Proof of Theorem \ref{thm}}]
	By Reineke's cell decomposition, we obtain that 
	$$
		{A_{S_*} = Z_{S_*} \cup \bigcup_{S''_* > S_*} Z_{S''_*}} = {\ZZ_{S_*} \cup \bigcup_{S''_* > S_*} \ZZ_{S''_*}}.
	$$ 
	The proper components of the intersection of the $D_{S_*}(x')$ with $x' \in C(S_*)$ are among those $\ZZ_{S''_*}$ with $S''_* \geq S_*$ and $d(S_*) = d(S''_*)$. Hence, using \cite[Ex.\ 8.2.1]{Fulton:98}, there are positive integers $n_{S_*,S_*}$ and $n_{S_*,S''_*}$ such that
	\begin{align*}
		\prod_{x' \in C(S_*)} c_{i(x')}(\UU) \cap [\Hilb_{d,n}]\ = \prod_{x' \in C(S_*)} [D_{S_*}(x')] 
		= n_{S_*,S_*} [\ZZ_{S_*}]\ + \sum_{\substack{S''_* > S_*,\\ d(S''_*) = d(S_*)}}\!\!\!\! n_{S_*,S''_*}[\ZZ_{S''_*}].
	\end{align*}
	It remains to prove that the coefficient $n_{S_*,S_*}$ is $1$. In order to do so, it suffices to prove that $D_{S_*} \cap U_{S_*} = Z_{S_*}$ as schemes. As mentioned above, Reineke shows that an isomorphism $U_{S*} \to \A^N$ (with $N := (m-1)d^2 + nd$) is given by the functions $\lambda_{x,x'}$ with $x \in S_*$ and $x' \in C(S_*)$ assigning to every point $[f,\phi]$ of $U_{S_*}$ the coefficient $\lambda_{x,x'}[f,\phi]$ that occurs displaying $\phi_{w'}fe_{k'}$ as a linear combination
	$$
		\phi_{w'}fe_{k'} = \sum_{x = (k,w) \in S_*} \lambda_{x,x'}[f,\phi] \cdot \phi_wfe_k
	$$
	where $x' = (k',w')$. Over $U_{S_*}$, the bundle $\UU$ trivializes. Moreover, for every pair $x_0 = (k_0,w_0)$, the sections $s_{x_0}$ of $\UU$ correspond to the sections of the trivial rank $d$-bundle on $U_{S_*}$ assigning to $[f,\phi]$ the matrix $(a_{x,x_0} \mid x \in S_*)$ of coefficients of the linear combination $\phi_{w_0}fe_{k_0} = \smash{\sum_{x = (k,w) \in S_*} a_{x,x_0} \phi_wfe_k}$. Therefore, enumerating $S_* = \{x_1 < \cdots < x_d\}$, the section $s_{x_i}$ restricted to $\A^N$ maps a matrix $\lambda$ to the $i$-th coordinate vector and $s_{x'}$ maps $\lambda$ to the vector $(\lambda_{x_1,x'},\ldots,\lambda_{x_d,x'})$. Under the isomorphism $U_{S_*} \to \A^N$, the degeneracy locus is therefore given by the vanishing of all $j(x')$-minors of the matrices
	$$
		\begin{pmatrix}
			1	&		&		& \lambda_{x_1,x'} \\
				& \ddots	&		& \vdots \\
				&		& 1 		& \lambda_{x_{j(x')},x'} \\
				&		& 		& \vdots \\
				&		& 		& \lambda_{x_d,x'}
		\end{pmatrix}.
	$$
	Thus, locally on $U_{S_*} \cong \A^N$, the degeneracy locus $D_{S_*}$ is given by the ideal generated by the coordinate functions $\smash{\lambda_{x_{j(x')+1},x'}},\ldots,\lambda_{x_d,x'}$ with $x'$ ranging over all critical pairs of $S_*$. It is therefore an affine space of dimension $d(S_*)$ and thus isomorphic to $Z_{S_*}$.
\end{proof}

\begin{rem} \label{rem_inter_mult}
	Let's fix the notation as in Theorem \ref{thm}. 
	\begin{enumerate} 
		\item We are able to determine the numbers $n_{S_*,S''_*}$ - at least in principle. They are given as intersection multiplicities
		$$
			n_{S_*,S''_*} = i \left( \ZZ_{S''_*}, D_{S_*}(x'_1) \cdots D_{S_*}(x'_r); \Hilb_{d,n}^{(m)} \right)
		$$
		as defined in \cite[Ex.\ 8.2.1]{Fulton:98}. Here, $\{x'_1,\ldots,x'_r\} = C(S_*)$. As the non-commutative Hilbert scheme is non-singular, it is also Cohen-Macaulay.  Applying Lemma \ref{lem_codim} and \cite[Ex.\ 14.4.2]{Fulton:98}, we obtain that every $D_{S_*}(x')$ is Cohen-Macaulay, too. Thus, \cite[Prop.\ 8.2, Ex.\ 8.2.7]{Fulton:98} imply that
		$$
			n_{S_*,S''_*} = l\left( \OO_{D_{S_*},\ZZ_{S''_*}} \right) = l\left(\OO_{D_{S_*} \cap U_{S''_*},Z_{S''_*}}\right),
		$$
		$D_{S_*}$ being equipped with its natural scheme structure.
		\item As $D_{S_*} = \smash{\bigcup_{x' \in C(S_*)}} D_{S_*}(x')$ and as every $D_{S_*}(x')$ has pure dimension $N-i(x')$, we see that every irreducible component of $D_{S_*}$ has dimension at least $N - \sum_{x'} i(x') = d(S_*)$ which is the dimension of $\ZZ_{S_*}$. The proper components of the intersection are precisely those with dimension $d(S_*)$. On the other hand, we know that 
		$$
			D_{S_*} = \bigsqcup_{S''_* \geq S_*} Z_{S''_*} = \bigcup_{S''_* \geq S_*} \ZZ_{S''_*}
		$$
		and that $\ZZ_{S''_*}$ are irreducible closed subsets of $D_{S_*}$ of dimension $d(S''_*) \leq d(S_*)$. For general reasons, the irreducible components of $D_{S_*}$ are those $\ZZ_{S''_*}$ (with $S''_* \geq S_*$) such that there exists no $S'_* \geq S_*$ with $\ZZ_{S''_*} \sub \ZZ_{S'_*}$ (other than $S'_* = S''_*$). We can easily see that these are exactly the sets $\ZZ_{S''_*}$ with $S''_* \geq S_*$ and $d(S''_*) = d(S_*)$. We have thus shown that the decomposition of $D_{S_*}$ into irreducible components reads as
		$$
			D_{S_*} = \bigcup_{\substack{S''_* > S_*,\\ d(S''_*) = d(S_*)}}\!\!\!\! \ZZ_{S''_*}
		$$
		which shows that $D_{S_*}$ has pure dimension $d(S_*)$.
		%
		\item Combining the previous two observations with Theorem \ref{thm}, we deduce that
		$$
			\prod_{x' \in C(S_*)} c_{i(x')}(\UU) \cap [\Hilb_{d,n}^{(m)}]\ =\ [D_{S_*}],
		$$
		the cycle attached to the closed subscheme $D_{S_*}$.
	\end{enumerate}
\end{rem}

\subsection{A basis consisting of monomials in Chern classes}

As a cell decomposition provides a basis of the Chow group of the variety, Theorem \ref{thm} gives a basis of the Chow group of the non-commutative Hilbert scheme consisting of monomials in the Chern classes of $\UU$. We want to know more concretely which monomials occur in this basis. Let $\smash{\mathcal{J}_{d,n}^{(m)}}$ be the set of tuples $j_* = (j_1,\ldots,j_{(m-1)d+n})$ of integers $0 \leq j_1 \leq \cdots \leq j_{(m-1)d+n} \leq d$ such that $j_\nu \geq l$ for every pair $(l,\nu)$ satisfying $1 \leq l \leq d$ and $(m-1)(l-1)+n \leq \nu$. If we display the non-decreasing sequence $j_*$ in a rectangular box, the condition for $j_*$ to belong to $\mathcal{J}_{d,n}$ precisely means that it has to cover the ``staircase'' as shown in the picture below:
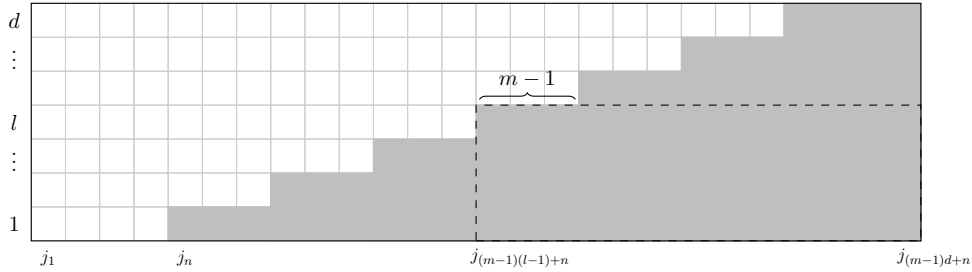
\begin{center}
	\begin{tikzpicture}[scale=.9]
	\foreach \x in {1,...,26}
	\foreach \y in {1,...,7}
	{
		\draw[gray!40, very thin] (\x*0.5,\y*0.5) +(-.25,-.25) rectangle ++(.25,.25);
	}
	\foreach \l in {1,...,7}
	{
		\filldraw[gray!50,fill=gray!50] (\l*1.5+1,0.5) +(-0.25,-0.25)  rectangle ++(12.25-\l*1.5,\l*0.5-0.25);
	}
	\draw (.25,.25) rectangle (13.25,3.75);
	\draw (0,0.5) node[scale=.75] {$1$};
	\draw (0,3*0.5) node[scale=.75] {$\vdots$};
	\draw (0,4*0.5) node[scale=.75] {$l$};
	\draw (0,6*0.5) node[scale=.75] {$\vdots$};
	\draw (0,7*0.5) node[scale=.75] {$d$};
	\draw (0.5,0) node[scale=.6] {$j_1$};
	\draw (5*0.5,0) node[scale=.6] {$j_n$};
	\draw (14*0.5+0.4,0) node[scale=.6] {$j_{(m-1)(l-1)+n}$};
	\draw (26*0.5+0.45,0) node[scale=.6] {$j_{(m-1)d+n}$};
	%
	\draw[dashed] (4*1.5+1,0.5) +(-0.25,-0.25)  rectangle ++(12.25-4*1.5,4*0.5-0.25);
	\draw decorate [decoration={name=brace}, yshift=.5ex]  {(14*0.5,5*0.5) +(-0.2,-0.25) -- node[above=0.4ex,scale=.75] {$m-1$} (16*0.5+0.2,5*0.5-0.25) };
\end{tikzpicture}
\end{center}
For a forest $S_* \in \FF_{d,n}$, we enumerate $S_* = \{x_1 < \cdots < x_d\}$ and denote $C(S_*) = \smash{\{x'_1, \ldots, x'_{(m-1)d+n}\}}$. To $S_*$ we associate the tuple $j_{S_*} = (j_1,\ldots,\allowbreak j_{(m-1)d+n})$ defined by letting $j_\nu = \smash{j_{S_*}(x'_{\nu})}$ be the maximal index $j$ for which $x_j < x'_{\nu}$.

\begin{lem} \label{lem_FJ}
	The mapping $S_* \mapsto j_{S_*}$ gives a bijection $\FF_{d,n}^{(m)} \to \mathcal{J}_{d,n}^{(m)}$.
\end{lem}

\begin{proof}
	We prove that the above map is well-defined and a bijection by induction on $d$, the number of nodes of the forest. For $d = 0$, the set $\FF_{d,n}$ consists of a single element, the empty forest, and $\mathcal{J}_{d,n}$ contains a single sequence, namely $(0,\ldots,0)$ of length $n$. Let $d > 0$. For a forest $S_* = \{ x_1 < \cdots < x_d \}$, we define $S'_* = \{ x_1 < \cdots < x_{d-1} \}$. We have
	$$
		C(S_*) = \left( C(S'_*) - \{x_d\} \right) \sqcup \{ \text{direct successors of } x_d \}.
	$$
	Enumerating $C(S_*) = \smash{\{x'_1 < \cdots < x'_{(m-1)d+n} \}}$ and letting $\nu_0$ be the maximal index $\nu$ such that $x'_\nu < x_d$ (we formally define $\nu_0 = 0$ if $x'_1 > x_d$), we see that $x'_1,\ldots,x'_{\nu_0}$ belong to $C(S'_*) - \{x_d\}$. This implies $\nu_0 \leq (m-1)(d-1)+n-1$. Then,
	$$
		j_\nu = j_{S_*}(x'_{\nu}) = \begin{cases} j_{S'_*}(x'_\nu) & \text{if } \nu \leq \nu_0 \\ d & \text{otherwise} \end{cases}
	$$
	for every $1 \leq \nu \leq (m-1)d+n$, which yields that the defining condition of $\mathcal{J}_{d,n}$ is fulfilled using that, by induction hypothesis, $j_{S'_*}$ belongs to $\mathcal{J}_{d-1,n}$. Next, we show that the map $\FF_{d,n} \to \mathcal{J}_{d,n}$ given by $S_* \mapsto j_{S_*}$ is bijective under the hypothesis that $\FF_{d-1,n} \to \mathcal{J}_{d-1,n}$ is. Let $j_* = (j_1,\ldots,j_{(m-1)d+n})$ be a sequence in $\mathcal{J}_{d,n}$. Define $\nu_0$ to be the maximal index $\nu$ such that $j_{\nu} < d$ (it will become clear a little later that this definition of $\nu_0$ agrees with the one before). Define $j'_* = (j'_1,\ldots,j'_{(m-1)(d-1)+n})$ by
	$$
		j'_\nu = \begin{cases} j_\nu & \text{if } \nu \leq \nu_0 \\ d-1 & \text{if } \nu > \nu_0. \end{cases}
	$$
	Then, $j'$ lies in $\mathcal{J}_{d-1,n}$, whence there exists a unique forest $S'_* \in \FF_{d-1,n}$ such that $j'_* = j_{S'_*}$. We order $C(S') = \smash{\{x''_1 < \cdots < x''_{(m-1)(d-1)+n} \}}$. If we construct $S_*$ by appending $x_d = x''_{\nu_0+1}$, we obtain a forest satisfying $j_{S_*} = j_*$. For any forest $S_*$ with $j_{S_*} = j_*$, we must have $j_{S'_*} = j'_*$. Therefore, as $S'_*$ is uniquely determined by $j'_*$, the given sequence $j_*$ determines $S_*$ uniquely.
\end{proof}

We introduce the set $\smash{\mathcal{B}_{d,n}^{(m)}}$ of all tuples $(b_0,\ldots,b_{d-1})$ of non-negative integers such that $\sum_{r=0}^i b_r < (m-1)i+n$ for all $0 \leq i \leq d-1$. For a sequence $j_* = (j_1,\ldots,j_{(m-1)d+n})$ in $\mathcal{J}_{d,n}$, we define $b_i$ to be the number of indexes $\nu$ with $j_\nu = i$. Denote $b_{j_*} = (b_0,\ldots,b_{d-1})$.

\begin{lem} \label{lem_JB}
	The association $j_* \mapsto b_{j_*}$ is a bijection $\mathcal{J}_{d,n}^{(m)} \to \mathcal{B}_{d,n}^{(m)}$.
\end{lem}

\begin{proof}
	Let's show that $b_{j_*}$ lies in $\mathcal{B}_{d,n}$. For every $i$, the sum $\sum_{r=0}^i b_r$ is the number of $\nu$ such that $j_{\nu} \leq i$. As, by definition of $\mathcal{J}_{d,n}$, all $j_{\nu}$ with $\nu \geq (m-1)i+n$ are greater or equal $i+1$, we obtain that $\sum_{r=0}^i b_r < (m-1)i+n$. It remains to prove that $\mathcal{J}_{d,n} \to \mathcal{B}_{d,n}$ is a bijection. We construct an inverse. Given a tuple $(b_0,\ldots,b_{d-1})$, we define $j_* = (j_1,\ldots,j_{(m-1)d+n})$ by
	$$
		j_{b_0 + \cdots + b_{i-1} + 1} = \cdots = j_{b_0 + \cdots + b_i} = i
	$$
	for all $0 \leq i \leq d$. Here, $b_d$ is defined as $(m-1)d+n - \sum_{r=0}^{d-1} b_r$. We see that $j_*$ belongs to $\mathcal{J}_{d,n}$ and that the maps $\mathcal{J}_{d,n} \to \mathcal{B}_{d,n}$ and $\mathcal{B}_{d,n} \to \mathcal{J}_{d,n}$ are mutually inverse.
\end{proof}

In combination with Theorem \ref{thm} and Corollary \ref{free}, the two previous lemmas prove the following result:

\begin{thm} \label{ko_basis_Chern}
	The Chow group $A_*(\Hilb_{d,n}^{(m)})$ is a free abelian group with basis
	$$
		c_1(\UU)^{b_{d-1}}\cdots c_d(\UU)^{b_0},
	$$
	where $(b_0,\ldots,b_{d-1})$ ranges over all tuples of non-negative integers satisfying $b_0 + \cdots + b_i < (m-1)i +n$ for every $0 \leq i \leq d-1$.
\end{thm}

In particular, the Chow \emph{ring} is generated by the Chern classes $c_1(\UU),\ldots,\allowbreak c_d(\UU)$. A result of King and Walter (cf.\ \cite[Thm.\ 3]{KW:95}) asserts that the Chow ring of a fine quiver moduli space is generated by the Chern classes of the universal bundle if the quiver has no oriented cycles. This theorem is not applicable here, yet the statement holds.

\begin{ex*}[continued]
	Let's illustrate Corollary \ref{ko_basis_Chern} in our favorite example $m=2$, $d=3$ and $n=1$. The bijections $\FF \to \mathcal{J} \to \mathcal{B}$ yield the following result: The trees
	
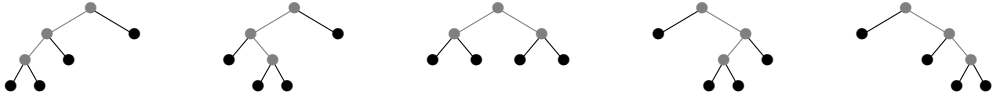
\begin{center}
		\begin{tikzpicture}[>=stealth',level/.style={sibling distance = 3em/#1, level distance = .9em}, parent anchor = center, child anchor = center] 
			\node (s) {}
				child { node (s1) {} edge from parent[gray]
					child { node (s11) {} edge from parent[gray]
						child { node (s111) {} edge from parent[black]}
						child { node (s112) {} edge from parent[black]}
					}
					child { node (s12) {} edge from parent[black]
						child { node (s121) {} edge from parent[white]}
						child { node (s122) {} edge from parent[white]}
					}
				}
				child { node (s2) {} edge from parent[black]
					child { node (s21) {} edge from parent[white]
						child { node (s211) {} edge from parent[white]}
						child { node (s212) {} edge from parent[white]}
					}
					child { node (s22) {} edge from parent[white]
						child { node (s221) {} edge from parent[white]}
						child { node (s222) {} edge from parent[white]}
					}
				}
			;
			\node [gray] at (s) {$\bullet$};
			\node [gray] at (s1) {$\bullet$};
			\node [gray] at (s11) {$\bullet$};
			\node [black] at (s111) {$\bullet$};
			\node [black] at (s112) {$\bullet$};
			\node [black] at (s12) {$\bullet$};
			\node [white] at (s121) {$\bullet$};
			\node [white] at (s122) {$\bullet$};
			\node [black] at (s2) {$\bullet$};
			\node [white] at (s21) {$\bullet$};
			\node [white] at (s211) {$\bullet$};
			\node [white] at (s212) {$\bullet$};
			\node [white] at (s22) {$\bullet$};
			\node [white] at (s221) {$\bullet$};
			\node [white] at (s222) {$\bullet$};
		\end{tikzpicture}
		\begin{tikzpicture}[>=stealth',level/.style={sibling distance = 3em/#1, level distance = .9em}, parent anchor = center, child anchor = center] 
			\node (s) {}
				child { node (s1) {} edge from parent[gray]
					child { node (s11) {} edge from parent[black]
						child { node (s111) {} edge from parent[white]}
						child { node (s112) {} edge from parent[white]}
					}
					child { node (s12) {} edge from parent[gray]
						child { node (s121) {} edge from parent[black]}
						child { node (s122) {} edge from parent[black]}
					}
				}
				child { node (s2) {} edge from parent[black]
					child { node (s21) {} edge from parent[white]
						child { node (s211) {} edge from parent[white]}
						child { node (s212) {} edge from parent[white]}
					}
					child { node (s22) {} edge from parent[white]
						child { node (s221) {} edge from parent[white]}
						child { node (s222) {} edge from parent[white]}
					}
				}
			;
			\node [gray] at (s) {$\bullet$};
			\node [gray] at (s1) {$\bullet$};
			\node [black] at (s11) {$\bullet$};
			\node [white] at (s111) {$\bullet$};
			\node [white] at (s112) {$\bullet$};
			\node [gray] at (s12) {$\bullet$};
			\node [black] at (s121) {$\bullet$};
			\node [black] at (s122) {$\bullet$};
			\node [black] at (s2) {$\bullet$};
			\node [white] at (s21) {$\bullet$};
			\node [white] at (s211) {$\bullet$};
			\node [white] at (s212) {$\bullet$};
			\node [white] at (s22) {$\bullet$};
			\node [white] at (s221) {$\bullet$};
			\node [white] at (s222) {$\bullet$};
		\end{tikzpicture}
		\begin{tikzpicture}[>=stealth',level/.style={sibling distance = 3em/#1, level distance = .9em}, parent anchor = center, child anchor = center] 
			\node (s) {}
				child { node (s1) {} edge from parent[gray]
					child { node (s11) {} edge from parent[black]
						child { node (s111) {} edge from parent[white]}
						child { node (s112) {} edge from parent[white]}
					}
					child { node (s12) {} edge from parent[black]
						child { node (s121) {} edge from parent[white]}
						child { node (s122) {} edge from parent[white]}
					}
				}
				child { node (s2) {} edge from parent[gray]
					child { node (s21) {} edge from parent[black]
						child { node (s211) {} edge from parent[white]}
						child { node (s212) {} edge from parent[white]}
					}
					child { node (s22) {} edge from parent[black]
						child { node (s221) {} edge from parent[white]}
						child { node (s222) {} edge from parent[white]}
					}
				}
			;
			\node [gray] at (s) {$\bullet$};
			\node [gray] at (s1) {$\bullet$};
			\node [black] at (s11) {$\bullet$};
			\node [white] at (s111) {$\bullet$};
			\node [white] at (s112) {$\bullet$};
			\node [black] at (s12) {$\bullet$};
			\node [white] at (s121) {$\bullet$};
			\node [white] at (s122) {$\bullet$};
			\node [gray] at (s2) {$\bullet$};
			\node [black] at (s21) {$\bullet$};
			\node [white] at (s211) {$\bullet$};
			\node [white] at (s212) {$\bullet$};
			\node [black] at (s22) {$\bullet$};
			\node [white] at (s221) {$\bullet$};
			\node [white] at (s222) {$\bullet$};
		\end{tikzpicture}
		\begin{tikzpicture}[>=stealth',level/.style={sibling distance = 3em/#1, level distance = .9em}, parent anchor = center, child anchor = center] 
			\node (s) {}
				child { node (s1) {} edge from parent[black]
					child { node (s11) {} edge from parent[white]
						child { node (s111) {} edge from parent[white]}
						child { node (s112) {} edge from parent[white]}
					}
					child { node (s12) {} edge from parent[white]
						child { node (s121) {} edge from parent[white]}
						child { node (s122) {} edge from parent[white]}
					}
				}
				child { node (s2) {} edge from parent[gray]
					child { node (s21) {} edge from parent[gray]
						child { node (s211) {} edge from parent[black]}
						child { node (s212) {} edge from parent[black]}
					}
					child { node (s22) {} edge from parent[black]
						child { node (s221) {} edge from parent[white]}
						child { node (s222) {} edge from parent[white]}
					}
				}
			;
			\node [gray] at (s) {$\bullet$};
			\node [black] at (s1) {$\bullet$};
			\node [white] at (s11) {$\bullet$};
			\node [white] at (s111) {$\bullet$};
			\node [white] at (s112) {$\bullet$};
			\node [white] at (s12) {$\bullet$};
			\node [white] at (s121) {$\bullet$};
			\node [white] at (s122) {$\bullet$};
			\node [gray] at (s2) {$\bullet$};
			\node [gray] at (s21) {$\bullet$};
			\node [black] at (s211) {$\bullet$};
			\node [black] at (s212) {$\bullet$};
			\node [black] at (s22) {$\bullet$};
			\node [white] at (s221) {$\bullet$};
			\node [white] at (s222) {$\bullet$};
		\end{tikzpicture}
		\begin{tikzpicture}[>=stealth',level/.style={sibling distance = 3em/#1, level distance = .9em}, parent anchor = center, child anchor = center] 
			\node (s) {}
				child { node (s1) {} edge from parent[black]
					child { node (s11) {} edge from parent[white]
						child { node (s111) {} edge from parent[white]}
						child { node (s112) {} edge from parent[white]}
					}
					child { node (s12) {} edge from parent[white]
						child { node (s121) {} edge from parent[white]}
						child { node (s122) {} edge from parent[white]}
					}
				}
				child { node (s2) {} edge from parent[gray]
					child { node (s21) {} edge from parent[black]
						child { node (s211) {} edge from parent[white]}
						child { node (s212) {} edge from parent[white]}
					}
					child { node (s22) {} edge from parent[gray]
						child { node (s221) {} edge from parent[black]}
						child { node (s222) {} edge from parent[black]}
					}
				}
			;
			\node [gray] at (s) {$\bullet$};
			\node [black] at (s1) {$\bullet$};
			\node [white] at (s11) {$\bullet$};
			\node [white] at (s111) {$\bullet$};
			\node [white] at (s112) {$\bullet$};
			\node [white] at (s12) {$\bullet$};
			\node [white] at (s121) {$\bullet$};
			\node [white] at (s122) {$\bullet$};
			\node [gray] at (s2) {$\bullet$};
			\node [black] at (s21) {$\bullet$};
			\node [white] at (s211) {$\bullet$};
			\node [white] at (s212) {$\bullet$};
			\node [gray] at (s22) {$\bullet$};
			\node [black] at (s221) {$\bullet$};
			\node [black] at (s222) {$\bullet$};
		\end{tikzpicture}
	\end{center}
	again displayed together with the sets $C(S_*)$, give rise to the following sequences of numbers in $\mathcal{J}_{3,1}^{(2)}$
	$$
		3333 \hspace{5em} 2333 \hspace{5em} 2233 \hspace{5em} 1333 \hspace{5em} 1233
	$$
	which, in turn, correspond to the following sequences of integers in $\BB_{3,1}^{(2)}$
	$$
		000 \hspace{5.5em} 001 \hspace{5.5em} 002 \hspace{5.5em} 010 \hspace{5.5em} 011.
	$$
	By forming the sequences of partial sums and then increasing every entry by one, we end up at item (s) in Stanley's list of interpretations of the Catalan numbers (cf.\ \cite[Ex.\ 6.19]{Stanley:99}). We have thus obtained a basis for the Chow group of the non-commutative Hilbert scheme. It reads
	$$
		A_*(\Hilb_{3,1}^{(2)}) = \Z \cdot 1 \oplus \Z \cdot c_1(\UU) \oplus \Z \cdot c_1(\UU)^2 \oplus \Z \cdot c_2(\UU) \oplus \Z \cdot c_1(\UU)c_2(\UU).
	$$
	Applying Theorem \ref{thm} also gives us a relation between these (monomials in) Chern classes and the cycles associated to the cell closures. The theorem tells us that
	\begin{align*}
		1 \cap [\Hilb_{3,1}^{(2)}] &= [\ZZ_{\treeA}],\\
		c_1(\UU) \cap [\Hilb_{3,1}^{(2)}] &= [\ZZ_{\treeB}],\\
		c_1(\UU)^2 \cap [\Hilb_{3,1}^{(2)}] &= [\ZZ_{\treeC}] + n[\ZZ_{\treeD}],\\
		c_2(\UU) \cap [\Hilb_{3,1}^{(2)}] &= [\ZZ_{\treeD}], \text{ and}\\
		c_1(\UU)c_2(\UU) \cap [\Hilb_{3,1}^{(2)}] &= [\ZZ_{\treeE}]
	\end{align*}
	for some positive integer $n$. As we have remarked (cf.\ Remark \ref{rem_inter_mult}), the integer $n$ is precisely the intersection multiplicity 
	$$
		i\left( \ZZ_{\treeD},D_{\treeC}(11) \cdot D_{\treeC}(12); \Hilb_{3,1}^{(2)} \right)
	$$
	which equals the length of the artinian local ring $\smash{\OO_{D_{\treeC} \cap U_{\treeD},Z_{\treeD}}}$. Employing the isomorphism $\smash{U_{\treeD}} \cong \A^{12}$ from \cite{Reineke:05}, we regard an element $[v,A,B] \in U_{\treeD}$ as a tuple\enlargethispage{1em}
	$$
		\begin{array}{c|ccc|ccc}
			1 & x_{\epsilon,1}	& 0 & x_{\epsilon,211}	& 0 & x_{\epsilon,22}	& x_{\epsilon,212} \\
			0 & x_{2,1}		& 1 & x_{2,211}		& 0 & x_{2,22}		& x_{2,212} \\
			0 & x_{21,1}		& 0 & x_{21,211}	& 1 & x_{21,22}		& x_{21,212}
		\end{array}
		=: %
		\begin{array}{c|ccc|ccc}
			1 & x_1		& 0 & x_2	& 0 & x_5	& x_8 \\
			0 & y		& 1 & x_3	& 0 & x_6	& x_9 \\
			0 & z		& 0 & x_4	& 1 & x_7	& x_{10}
		\end{array}
	$$
	by displaying $[v,A,B]$ in terms of the basis $v,Bv,ABv$. The closed subscheme $D_{\treeC} \cap U_{\treeD}$ is defined by the vanishing of the determinants
	\begin{align*}
		\det( v \mid Av \mid A^2v) & = %
		\begin{vmatrix}
			1 & x_1 & x_1^2 + x_2z \\
			0 & y   & x_1y + x_3z \\
			0 & z   & x_1z + y + x_4z \\
		\end{vmatrix}
		= y^2 + x_4yz - x_3z^2 \text{ and} \\
		\det( v \mid Av \mid BAv) &= %
		\begin{vmatrix}
			1 & x_1 & x_5y + x_8z \\
			0 & y   & x_1 + x_6y + x_{9}z \\
			0 & z   & x_7y + x_{10}z \\
		\end{vmatrix}\\
		&= x_7y^2 + (x_{10}-x_6)yz - x_1z - x_{9}z^2.
	\end{align*}
	On the other hand, the closed subvariety $\smash{Z_{\treeD}}$ is given by the vanishing of $y$ and $z$. Therefore, the local ring of $D_{\treeC} \cap U_{\treeD}$ along the closed subvariety $Z_{\treeD}$ is
	$$
		\kk(x_1,\ldots,x_{10})[y,z]/(y^2 + x_4yz - x_3z^2,\ x_7y^2 + (x_{10}-x_6)yz - x_1z - x_{9}z^2).
	$$
	A lengthy computation shows that the length of this (artinian) ring - which equals its dimension over $\kk(x_1,\ldots,x_{10})$ - is $4$. The author has determined this using \textsc{Singular}.
\end{ex*}

	\section{A module structure over the cohomological Hall algebra}

The Cohomological Hall algebra, which we will call CoHa in the following, was invented by Kontsevich and Soibelman in \cite{KS:11}. We will consider the CoHa for the $m$-loop quiver and define a module structure on the Chow rings of non-commutative Hilbert schemes. In the first version of this article, we did this in the case $n=1$, but the results and proofs for arbitrary $n$ are completely analogous.

\subsection{Cohomological Hall algebra of a loop quiver}

We consider the $m$-loop quiver $Q$. A dimension vector for $Q$ is just a non-negative integer $d$. Denote by $R_d$ the vector space $M_{d \times d}^m$ of $m$-tuples of $(d \times d)$-matrices. On this space, we have an action of the reductive linear algebraic group $G_d := \Gl_d$ by conjugation. We define
$$
	\HH_d := A^*_{G_d}(R_d)_{\Q}
$$
the $G_d$-equivariant Chow ring as defined by Edidin and Graham (cf.\ \cite{EG:98}). See also Brion's article \cite{Brion:98}. Note that, in case $\kk = \C$, the equivariant cycle map $A^*_{G_d}(R_d)_{\Q} \to H^*_{G_d}(R_d;\Q)$ is an isomorphism (of rings) that doubles the degrees, as $R_d$ is a vector space with a linear $G_d$-action, whence $A_{G_d}^ *(R_d)_{\Q} \cong A_{G_d}^*(\pt)_{\Q}$ and $H_{G_d}^*(R_d;\Q) \cong H_{G_d}^*(\pt;\Q)$ and moreover, both $A_{G_d}(\pt)_{\Q}$ and $H_{G_d}(\pt;\Q)$ are isomorphic to the ring of symmetric polynomials $\Q[x_1,\ldots,x_d]^{S_d}$ in $d$ variables (cf.\ \cite[Subsect.\ 3.2]{EG:98}). Under these identifications, the cycle map corresponds to the identity. Therefore, we may work with equivariant Chow rings (with rational coefficients) rather than with equivariant cohomology rings (like Kontsevich--Soibelman do). We make the following convention: In this section, all (equivariant) Chow rings are extended to the rationals, i.e. write $A_G^*(X)$ instead of $A_G^*(X)_\Q$. On the direct sum $\HH := \bigoplus_{d \geq 0} \HH_d$, Kontsevich--Soibelman define a ``convolution like'' multiplication $\HH_p \otimes \HH_q \to \HH_{d = p+q}$, assigning $f \otimes g \mapsto f * g$, as the composition of the horizontal maps in
\begin{center}
	\begin{tikzpicture}[description/.style={fill=white,inner sep=2pt}]
		\matrix(m)[matrix of math nodes, row sep=2em, column sep=1.5em, text height=1.5ex, text depth=0.25ex]
		{
			A_{G_p}^i(R_p) \otimes A_{G_q}^j(R_q) \\
			A_{G_p \times G_q}^{k= i+j}(R_p \times R_q) &[-1.5em] \cong &[-1.5em] A_{P_{p,q}}^k(R_{\pqstar}) & A_{P_{p,q}}^{k+mpq}(R_d) & A_{G_d}^{k+(m-1)pq}(R_d). \\
			A_{G_p \times G_q}^k(R_{\pqstar}) & & & \\
		};
		\path[->, font=\scriptsize]
		(m-1-1) edge node[auto] {$\times$} (m-2-1)
		(m-2-1) edge node[auto] {$\cong$} (m-3-1)
		(m-2-3) edge node[auto] {$\cong$} (m-3-1)
		(m-2-3) edge node[auto] {$i_*$} (m-2-4)
		(m-2-4) edge node[auto] {$\pi_*$} (m-2-5)
		;
	\end{tikzpicture}
\end{center}
Here, $P_{p,q}$ denotes the upper parabolic of $G_d$ to the decomposition of $\kk^d$ into the coordinate space of the first $p$ and the last $q$ unit vectors. The subspace $R_{\pqstar}$ of $R_d$ is the space of all tuples $(\phi_1,\ldots,\phi_m)$ of $(d \times d)$-matrices, mapping the coordinate space $\kk^p$ into itself. From now on, write $L_{p,q} := G_p \times G_q$. It is the Levi factor of the parabolic $P_{p,q}$. The above maps arise as follows: The map $\times$ is the equivariant exterior product. For the other maps, choose a $G_d$-module $V$ such that there exists an open $G_p$-equivariant subset $U \sub V$ on which a fiber bundle quotient $U \to U/G_d$ exists and with the property that $\dim V - \dim(V-U) > k + mpq$. We have morphisms
\begin{center}
	\begin{tikzpicture}[description/.style={fill=white,inner sep=2pt}]
		\matrix(m)[matrix of math nodes, row sep=2em, column sep=1.7em, text height=1.5ex, text depth=0.25ex]
		{
			(R_p \times R_q) \times^{L_{p,q}} U &[-.7em]  &[-1em] R_{\pqstar} \times^{P_{p,q}} U & R_d \times^{P_{p,q}} U & R_d \times^{G_d} U, \\
			R_{\pqstar} \times^{L_{p,q}} U & & & \\
		};
		\path[->, font=\scriptsize]
		(m-2-1) edge (m-1-1)
		(m-2-1) edge (m-1-3)
		(m-1-3) edge node[auto] {$i$} (m-1-4)
		(m-1-4) edge node[auto] {$\pi$} (m-1-5)
		;
	\end{tikzpicture}
\end{center}
the non-horizontal maps being affine space bundles, $i$ is a closed embedding and $\pi$ is a smooth morphism with fiber $G_d/P_{p,q} = \Gr_p(\kk^d)$. In particular, $\pi$ is proper. Note that this choice of $U$ assures that the (ordinary) Chow groups of the varieties in the above diagram precisely yield the equivariant Chow groups occurring in the diagram before. The thus obtained multiplication makes $\HH$ an associative graded algebra with a unit $1 \in \HH_0$. Moreover, we can define a $(\Z_{\geq 0} \times \Z)$-bigrading on $\HH$ by putting
$$
	\HH_{d,k} := A^{(k-d^2)/2}_{G_d}(R_d)
$$
if $k - d^2$ is even and zero otherwise. This bigrading is coincides with the one in \cite[Subsect.\ 2.6]{KS:11}.

\begin{defn*}[{\cite[Def.\ 1]{KS:11}}]
	The bigraded algebra $\HH = \bigoplus_{d,k} \HH_{d,k}$ is called the \textbf{Cohomological Hall algebra} of the $m$-loop quiver.
\end{defn*}

In \cite[Thm.\ 2]{KS:11}, an explicit formula for the multiplication is derived. Identifying $\HH_d = A^*_{G_d}(R_d)$ with $\smash{A^*_{G_d}(\pt) \cong A^*_{T_d}(\pt)^W = \Q[x_1,\ldots,x_d]^{S_d}}$, we obtain
\begin{align*}
	(f*g)(x_1,\ldots,x_d) = \sum_{\substack{1 \leq i_1 < \cdots < i_p \leq d\\1 \leq j_1 < \cdots < j_q \leq d\\\text{complementary}}} f(x_{i_1},\ldots,x_{i_p})g(x_{j_1},\ldots,x_{j_q}) \prod_{\mu = 1}^p \prod_{\nu = 1}^q (x_{j_\nu} - x_{i_\mu})^{m-1},
\end{align*}
where two sequences $1 \leq i_1 < \cdots < i_p \leq d$ and $1 \leq j_1 < \cdots < j_q \leq d$ are called complementary if the union of these numbers is $\{1,\ldots,d\}$. Using this formula, it is evident that the multiplication $*$ is graded commutative if $m$ is even, and commutative if $m$ is odd.

\subsection{A computational approach}

Our goal is to realize $\bigoplus_d A^*(\Hilb_{d,n}^{(m)})$ as an $\HH$-module. In fact, it will turn out to be a quotient of $\HH$ by some ideal, thus it inherits the structure of an algebra itself. The idea stems from a purely algebraic observation. A point $(f,\phi) \in \hat{R}_{d,n}$ is stable if and only if the image of $f$ generates $\kk^d$ as a $\kk\langle \phi_1,\ldots,\phi_m \rangle$-module. In other words, this means that a proper subspace $U$ containing $\im f$ can't be invariant under all $\phi_i$. Consider the universal bundle $\UU$ on $\Hilb_{d,n} := \smash{\Hilb_{d,n}^{(m)}}$ of rank $d$. Let $\Fl := \Fl(\UU) \to \Hilb_{d,n}$ be the complete flag bundle. It possesses a universal flag
$$
	\UU^\bull: 0 = \UU^0 \sub \UU^1 \sub \cdots \sub \UU^d = \UU_{\Fl}
$$
with $\rk \UU^i = i$. A point $y$ of $\Fl(\UU)$ is a pair consisting of $[f,\phi] \in \Hilb_{d,n}$ and a complete flag $W^\bull$ of $\UU_{[f,\phi]} \cong \kk^d$. The universal flag is defined by the property that its fiber $\UU^p_{y}$ in the point $y$ is precisely $W^p$ for all points $y = ([f,\phi],W^\bull)$ of $\Fl$. Fix an integer $0 \leq p < d$. The section $s_k = s_{k,\epsilon}$ of $\UU$ readily induces a section $s_k^p$ of $\UU_{\Fl}/\UU^p$. Moreover, the endomorphisms $\Phi_1,\ldots,\Phi_m$ give homomorphisms
$$
	\UU^p \into \UU_{\Fl} \xto{}{\Phi_i} \UU_{\Fl} \to \UU_{\Fl}/\UU^p
$$
which may be interpreted as sections $\Phi_i^p$ of the Hom-bundle $\underline{\Hom}(\UU^p, \UU_{\Fl}/\UU^p) = \UU_{\Fl}/\UU^p \otimes (\UU^p)^\vee$. For all $y = ([f,\phi],W^\bull) \in \Fl$ as before, the subspace $W^p$ can't be invariant under all $\phi_i$ if it contains $\im f$. In terms of the above sections, this means that
$$
	s_1^p(y),\ldots,s_n^p(y),\Phi_1^p(y),\ldots,\Phi_m^p(y)
$$
can't all be zero. In other words, the intersection $Z(s_1^p) \cap \cdots \cap Z(s_n^p) \cap Z(\Phi_1^p) \cap \cdots \cap Z(\Phi_m^p)$ of the zero loci of the sections is empty. In particular, this implies that
$$
	\Z(s_1^p) \cdots \Z(s_n^p) \cdot \Z(\Phi_1^p)\cdots\Z(\Phi_m^p) = 0,
$$
$\Z$ denoting the localized top Chern class as in \cite[14.1]{Fulton:98}. The image of the left-hand expression in the Chow ring $A^*(\Fl)$ equals
\begin{align*}
	c_{\mathrm{top}}\left( \UU_{\Fl}/\UU^p \right)^n \cdot c_{\mathrm{top}}\left( \UU_{\Fl}/\UU^p \otimes (\UU^p)^\vee \right)^m &= \xi_{p+1}^n\cdots \xi_d^n \prod_{\mu=1}^p \prod_{\nu = p+1}^d (\xi_\nu - \xi_\mu)^m \\ 
	&=: f^{(p)}(\xi_1,\ldots,\xi_d)
\end{align*}
if we denote $\xi_\nu := c_1(\UU^\nu/\UU^{\nu-1})$. Like in \cite[Def.\ 12]{Franzen:15:Chow_Ring_Quiv:Publ}, we call $f^{(p)}$ a forbidden polynomial. By a result of Grothendieck (cf.\ \cite[Thm.\ 1]{Groth:58:inter}), the Chow ring $A^*(\Fl)$ is isomorphic to $A^*(\Hilb_{d,n}) \otimes_{\HH_d} C_d$ as a $C_d$-algebra, where $C_d = \Q[x_1,\ldots,x_d]$ is a polynomial ring and $\HH_d = \Q[x_1,\ldots,x_d]^{S_d} = \Q[e_1,\ldots,e_d]$ is the ring of symmetric polynomials. It is a basic fact that $C_d$ is a free $\HH_d$-module. Pick a basis $\BB$. Displaying $f^{(p)}$ in terms of $\BB$ we obtain coefficients $\tau^{(p)}(y)$ for every $y \in \BB$. As $f^{(p)}$ vanishes at the Chern \emph{roots} of $\UU$, the coefficient $\tau^{(p)}(y)$ (which is a polynomial in the elementary symmetric functions $e_1,\ldots,e_d$) vanishes at the Chern \emph{classes} of $\UU$. In \cite[Def.\ 13]{Franzen:15:Chow_Ring_Quiv:Publ}, $\tau^{(p)}(y)$ is called a tautological relation. We might wonder if these tautological relations provide a presentation of the Chow ring of $\Hilb_{d,n}$.
Let $C_d := \Q[x_1,\ldots,x_d]$ and let $\rho = \rho_d: C_d \to \HH_d = \Q[x_1,\ldots,x_d]^{S_d}$ be the $\HH_d$-linear map defined by
$$
	\rho(f) = \Delta^{-1} \sum_{w \in S_d} \sign(w) w.f.
$$
In this context, $\Delta = \Delta_d$ is the discriminant $\prod_{i < j} (x_j - x_i)$. In \cite{Franzen:15:Chow_Ring_Quiv:Publ}, it is shown that the ideal $(\tau^{(p)}(y) \mid y)$ generated by the tautological relations coincides with $\smash{\rho(C_d \cdot f^{(p)})}$. Thus, generators of this ideal are given by the images $\rho(b \cdot f^{(p)})$, where $b$ runs through a basis of $C_d$ over $\HH_d$. We now use the fact that for every $w \in S_d$, there exists a unique $\tau \in S_p \times S_q$ and a unique $(p,q)$-shuffle permutation $\sigma$ with $w = \sigma \tau$. Being a $(p,q)$-shuffle permutation means that $\sigma$ is of the form
$$
	\sigma = \begin{pmatrix}
	         	1 & \cdots & p & p+1 & \cdots & d \\
	         	i_1 & \cdots & i_p & j_1 & \cdots & j_q
	         \end{pmatrix}
$$
for sequences $i_1 < \cdots < i_p$ and $j_1 < \cdots < j_q$ which are necessarily complementary. For any $b \in C$, we obtain
\begin{align*}
	\rho(bf^{(p)}) &= \Delta^{-1} \sum_{\sigma\ (p,q)\text{-shuffle}}\ \sum_{\tau \in S_p \times S_q} \sign(\sigma \tau) (\sigma \tau).(bf^{(p)}) \\
	&= \Delta^{-1} \sum_{\sigma\ (p,q)\text{-shuffle}} \sign(\sigma) \sigma.f^{(p)} \sum_{\tau \in S_p \times S_q} \sign(\tau) (\sigma \tau).b \\
	&= \Delta^{-1} \sum_{\sigma\ (p,q)\text{-shuffle}} \sign(\sigma) \sigma.f^{(p)} \cdot \sigma.\!\!\left( \sum_{\tau \in S_p \times S_q} \sign(\tau) \tau.b \right).
\end{align*}
As $\sum_\tau \sign(\tau) \tau.b$ is alternating under the action of $S_p \times S_q$, it is divisible by $\Delta_{p \times q}$ which we define to be $\Delta_p(x_1,\ldots,x_p) \Delta_q(x_{p+1},\ldots,x_d)$. Setting $\delta$ as the product $\prod_{i = 1}^p \prod_{j = p+1}^d (x_j - x_i)$, we obtain $\Delta = \delta \Delta_{p \times q}$, and therefore,
$$
	\sign(\sigma) \Delta = \sigma.\Delta = (\sigma.\delta)(\sigma.\Delta_{p \times q}).
$$
With $\rho_{p \times q}(b) = \Delta_{p \times q}^{-1} \sum_{\tau \in S_p \times S_q} \sign(\tau) \tau.b$, we get
$$
	\rho(bf^{(p)}) = \sum_{\sigma\ (p,q)\text{-shuffle}} (\sigma.\delta)^{-1} \cdot \sigma.f^{(p)} \cdot \sigma.\rho_{p \times q}(b).
$$
We insert the definition of $f^{(p)}$. We write it as $f^{(p)} = x_{p+1}^n\cdots x_d^n \delta^m$. Thus, $\rho(bf^{(p)})$ equals
\begin{align*}
	\sum_{\substack{1 \leq i_1 < \cdots < i_p \leq d\\1 \leq j_1 < \cdots < j_q \leq d\\\text{complementary}}} x_{j_1}^n\cdots x_{j_q}^n \prod_{\mu = 1}^p \prod_{\nu = 1}^q &(x_{j_\nu} - x_{i_\mu})^{m-1} \cdot (\rho_{p \times q} (b))(x_{i_1},\ldots, x_{i_p},x_{j_1},\ldots,x_{j_q}).
\end{align*}
As $b$ runs through $C_d$, the image $\rho_{p \times q}(b)$ runs through $\Q[x_1,\ldots,x_p,x_{p+1},\ldots,\allowbreak x_d]^{S_p \times S_q}$ which we may identify with the tensor product $\Q[x_1,\ldots,x_p]^{S_p} \otimes \Q[x_{p+1},\ldots,x_d]^{S_q}$. Thus: 

\begin{lem}
The ideal $\rho(C_d\cdot f^{(p)})$ consists of expressions
	$$
		\sum_{\substack{1 \leq i_1 < \cdots < i_p \leq d\\1 \leq j_1 < \cdots < j_q \leq d\\\text{complementary}}} f(x_{i_1},\ldots,x_{i_p})g(x_{j_1},\ldots,x_{j_q})x_{j_1}^n\cdots x_{j_q}^n \prod_{\mu = 1}^p \prod_{\nu = 1}^q (x_{j_\nu} - x_{i_\mu})^{m-1},
	$$
	with $f$ ranging over all symmetric polynomials in $p$ variables and $g$ over those in $q$ variables.
\end{lem}

The similarity of the above term with the multiplication in the CoHa is too obvious to be coincidental. By abuse of notation, we write $e_i \in \HH_q$ for the $i$-th elementary symmetric function in $q$ variables whenever it is obvious where these elements live. We obtain that $\rho(C_d \cdot f^{(p)})$ is generated as a rational vector space by all expressions $f * (e_q^n \cup g)$ with $f \in \HH_p$ and $g \in \HH_q$, writing $\cup$ for the ordinary multiplication in the ring of symmetric functions (which coincides with the cup product) in order to distinguish it from the CoHa-multiplication. Note that $e_q^n$ is the $n$-th power of $e_q$ with respect to the cup product. Summarizing, we have seen:

\begin{lem} \label{lem_taut}
\hspace{-1ex}The ideal $\sum_{p < d} \sum_{y \in \BB} \HH_d \cdot \tau^{(p)}(y) = \sum_{p < d} \rho(C_d \cdot f^{(p)})$ of tautological relations equals
	$$
		\sum_{p+q = d,\ q \neq 0} \HH_p * (e_q^n \cup \HH_q).
	$$
\end{lem}

\subsection{Construction of the CoHa-module structure}

Put $\AA := \bigoplus_d \AA_d$, where $\AA_d := A^*(\Hilb_{d,n})$ is the Chow ring of the non-commutative Hilbert scheme of codimension $d$ ideals. We are going to define an $\HH$-module structure on $\AA$ and show that we can realize it as a quotient of $\HH$ itself and describe the kernel of the quotient map. Note that we may construct a similar diagram as above for the ``decorated'' representation variety $\hat{R}_{d,n} = M_{d \times n} \times M_{d \times d}^m$ (cf.\ also page \pageref{page_def_Rhat}): We have morphisms
\begin{center}
	\begin{tikzpicture}[description/.style={fill=white,inner sep=2pt}]
		\matrix(m)[matrix of math nodes, row sep=2em, column sep=1.7em, text height=1.5ex, text depth=0.25ex]
		{
			(R_p \times \hat{R}_{q,n}) \times^{L_{p,q}} U &[-.7em]  &[-1em] \hat{R}_{\pqstar,n} \times^{P_{p,q}} U & \hat{R}_{d,n} \times^{P_{p,q}} U & \hat{R}_{d,n} \times^{G_d} U, \\
			\hat{R}_{\pqstar,n} \times^{L_{p,q}} U & & & \\
		};
		\path[->, font=\scriptsize]
		(m-2-1) edge (m-1-1)
		(m-2-1) edge (m-1-3)
		(m-1-3) edge node[auto] {$i$} (m-1-4)
		(m-1-4) edge node[auto] {$\pi$} (m-1-5)
		;
	\end{tikzpicture}
\end{center}
where the arrows without names are again affine space bundles. This induces an $\HH$-module structure on the direct sum $\bigoplus_d A^*_{G_d}(\hat{R}_{d,n})$. But as $\hat{R}_{d,n} \to R_d$ is also a $G_d$-equivariant affine space bundle, the direct sum of these Chow groups coincides with $\HH$, as a vector space. It is not hard to see that the induced module structure on $\smash{\bigoplus_d A^*_{G_d}(\hat{R}_{d,n})}$ coincides with the natural $\HH$-module structure on $\HH$ itself. As a next step, we pass to the stable locus of $\smash{\hat{R}_{d,n}}$. It consists of the tuples $(f,\phi)$ such that $\im f$ generates $\kk^d$ as a left-$\kk\langle \phi_1,\ldots,\phi_m \rangle$-module. Consider the open subset 
$$
	\hat{R}_{\pqstar,n}^{\st} = \hat{R}_{\pqstar,n} \cap \hat{R}_{d,n}^{\st}
$$
of $\smash{\hat{R}_{\pqstar,n} = M_{d \times n} \times R_{\pqstar}}$. An element $(f,\phi)$ of $\smash{\hat{R}_{\pqstar,n}^{\st}}$ is of the form $\Big( \big(\begin{smallmatrix} f' \\ \hline f'' \end{smallmatrix}\big), \big(\begin{smallmatrix} \phi' & \vline & * \\ \hline 0 & \vline & \phi'' \end{smallmatrix}\big) \Big)$
and thus, for any $P \in \kk\langle x_1,\ldots,x_m \rangle$ and any $1 \leq k \leq n$, we obtain
\begin{align*}
	P(\phi_1,\ldots,\phi_m)fe_k &= 	\left( \begin{array}{c | c}
	                                                 	P(\phi'_1,\ldots,\phi'_m) & * \\
	                                                 	\hline
	                                                 	0 & P(\phi''_1,\ldots,\phi''_m)
	                              \end{array} \right)  
	                              \cdot \left(\begin{array}{c} f'e_k \\ \hline  f''e_k \end{array}\right) \\
	                              &= 
	                                                 \left(\begin{array}{c} * \\ \hline  P(\phi''_1,\ldots,\phi''_m)f''e_k \end{array}\right).
\end{align*}
This implies that $(f'',\phi'')$ is a stable point of $\hat{R}_{q,n}$. By restricting the projection $\hat{R}_{\pqstar,n} \to R_p \times \hat{R}_{q,n}$, we obtain a well defined morphism
$$
	F: \hat{R}_{\pqstar,n}^{\st} \to R_p \times \hat{R}_{q,n}^{\st}.
$$
Being an affine space bundle, the projection $\hat{R}_{\pqstar,n} \to R_p \times \hat{R}_{q,n}$ is flat. This implies at once that $F$ is a flat morphism, too. We can draw the diagram
\begin{center}
	\begin{tikzpicture}[description/.style={fill=white,inner sep=2pt}]
		\matrix(m)[matrix of math nodes, row sep=2em, column sep=1.7em, text height=1.5ex, text depth=0.25ex]
		{
			(R_p \times \hat{R}_{q,n}^{\st}) \times^{L_{p,q}} U &[-.7em]  &[-1em] \hat{R}_{\pqstar,n}^{\st} \times^{P_{p,q}} U & \hat{R}_{d,n}^{\st} \times^{P_{p,q}} U & \hat{R}_{d,n}^{\st} \times^{G_d} U, \\
			\hat{R}_{\pqstar,n}^{\st} \times^{L_{p,q}} U & & & \\
		};
		\path[->, font=\scriptsize]
		(m-2-1) edge node[auto] {$F$} (m-1-1)
		(m-2-1) edge (m-1-3)
		(m-1-3) edge node[auto] {$i$} (m-1-4)
		(m-1-4) edge node[auto] {$\pi$} (m-1-5)
		;
	\end{tikzpicture}
\end{center}
which gives maps in equivariant intersection theory as follows:
\begin{center}
	\begin{tikzpicture}[description/.style={fill=white,inner sep=2pt}]
		\matrix(m)[matrix of math nodes, row sep=2em, column sep=1em, text height=1.5ex, text depth=0.25ex]
		{
			A_{G_p}^i(R_p) \otimes A_{G_q}^j(\hat{R}_{q,n}^{\st})\\ 
			A_{G_p \times G_q}^{k}(R_p \times \hat{R}_{q,n}^{\st}) &[-.7em] &[-1em] A_{P_{p,q}}^k(\hat{R}_{\pqstar,n}^{\st}) & A_{P_{p,q}}^{k+mpq}(\hat{R}_{d,n}^{\st}) & A_{G_d}^{k+(m-1)pq}(\hat{R}_{d,n}^{\st}). \\
			A_{G_p \times G_q}^k(\hat{R}_{\pqstar,n}^{\st}) & & & \\
		};
		\path[->, font=\scriptsize]
		(m-1-1) edge node[auto] {$\times$} (m-2-1)
		(m-2-1) edge node[auto] {$F^*$} (m-3-1)
		(m-2-1) edge (m-2-3)
		(m-2-3) edge node[auto] {$\cong$} (m-3-1)
		(m-2-3) edge node[auto] {$i_*$} (m-2-4)
		(m-2-4) edge node[auto] {$\pi_*$} (m-2-5)
		;
	\end{tikzpicture}
\end{center}
Composing these maps, we get $\HH_p \otimes \AA_q \to \AA_d$. Let's write $f \circledast g$ for the image of $f \otimes g$ under this map. A similar argument to \cite[2.3]{KS:11} shows:

\begin{prop}
	The maps $\HH_p \otimes \AA_q \to \AA_{p+q}$ constructed above make $\AA$ into an $\HH$-module.
\end{prop}

\begin{rem*}
	In a recent article, Soibelman has pointed this out as well (cf.\ \cite{Soibelman:14}). In fact, he considers CoHa-modules arising from stable framed objects in a much more general context than we do here.
\end{rem*}

If we define a bigrading on $\AA$ by letting
$
	\AA_{d,k} := A_{G_d}^{(k-d^2)/2}(\hat{R}_{d,n}^{\st}),
$
we obtain that $\AA$ also becomes a bigraded $\HH$-module.

Let's look at the map $j^*: \HH \to \AA$ which is induced by the open embeddings $\hat{R}_{d,n}^{\st} \into \hat{R}_{d,n}$. It is clearly surjective. It is also $\HH$-linear as the following commutative diagram asserts:
\begin{center}
	\begin{tikzpicture}[description/.style={fill=white,inner sep=2pt}]
		\matrix(m)[matrix of math nodes, row sep=2.0em, column sep=1.1em, text height=1.5ex, text depth=0.25ex]
		{
				&[-1.5em] (R_p \times \hat{R}_{q,n}^{\st}) \times^{L_{p,q}} U &[-2em]	&[1.5em] \hat{R}_{\pqstar,n}^{\st} \times^{L_{p,q}} U \\
			(R_p \times \hat{R}_{q,n}) \times^{L_{p,q}} U &	& \hat{R}_{\pqstar,n} \times^{L_{p,q}} U	& \hat{R}_{\pqstar,n}^{\st} \times^{P_{p,q}} U \\
				&	& \hat{R}_{\pqstar,n} \times^{P_{p,q}} U	& \hat{R}_{d,n}^{\st} \times^{P_{p,q}} U \\
				&	& \hat{R}_{d,n} \times^{P_{p,q}} U		& \hat{R}_{d,n}^{\st} \times^{G_d} U \\
				&	& \hat{R}_{d,n} \times^{G_d} U.		& \\
		};
		\path[->, font=\scriptsize]
		(m-1-2) edge (m-2-1)
		(m-1-4) edge node[above] {$F$} (m-1-2)
		(m-1-4) edge (m-2-3)
		(m-2-3) edge (m-2-1)
		(m-1-4) edge (m-2-4)
		(m-2-4) edge node[auto] {$i$} (m-3-4)
		(m-3-4) edge node[auto] {$\pi$} (m-4-4)
		(m-2-3) edge (m-3-3)
		(m-3-3) edge node[auto] {$i$} (m-4-3)
		(m-4-3) edge node[auto] {$\pi$} (m-5-3)
		(m-2-4) edge (m-3-3)
		(m-3-4) edge (m-4-3)
		(m-4-4) edge (m-5-3)
		;
	\end{tikzpicture}
\end{center}
In this diagram, all maps pointing from north-east to south-west are induced by the open embeddings. Note that every ``square'', except for the uppermost, is cartesian. Moreover, the upper ``square'' is commutative as $F$ is defined to be the restriction of the projection $\smash{\hat{R}_{\pqstar,n} \to R_p \times \hat{R}_{q,n}}$. Hence, passing to intersection theory, the outer arrows of the diagram give two ways to go from $\smash{A^*_{L_{p,q}}(R_p \times \hat{R}_{q,n})}$ to $\smash{A^*_{G_d}(\hat{R}_{d,n}^{\st})}$ which coincide. One way describes $f \circledast j^*g$, whereas the other represents $j^*(f * g)$. In a picture:
\begin{center}
	\begin{tikzpicture}
		\draw (0,2) node(a) {$A^*_{L_{p,q}}(R_p \times \hat{R}_{q,n})$} (5,1) node(b) {$A^*_{G_d}(\hat{R}_{d,n}^{\st})$.};
		\draw[->] (a.north east) [rounded corners] -- (3,3) -- (5,3) -- (b.north);
		\draw[->] (a.east) [rounded corners] -- (3,2) -- (3,0) -- (b.south west);
		\draw (6,2.5) node(c) {$f \circledast j^*g$};
		\draw (2,1) node(c) {$j^*(f * g)$};
	\end{tikzpicture}
\end{center}
Considering the above defined bigrading on $\AA$, the map $j^*:\HH \to \AA$ is also homogeneous of bidegree $(0,0)$. Summarizing:

\begin{prop}
	The map $j^*: \HH \to \AA$ induced by the open embeddings $j: \smash{\hat{R}_{d,n}^{\st}} \to \smash{\hat{R}_{d,n}}$ is $\HH$-linear, surjective and homogeneous of bidegree $(0,0)$.
\end{prop}

In other words, $\AA$ can be written as a quotient of $\HH$. Taking into account that $\HH$ is either commutative (if $m$ is odd) or graded commutative (if $m$ is even), we obtain:

\begin{cor}
	The vector space $\AA$ inherits the structure of a bigraded $\HH$-algebra.
\end{cor}

Motivated by the calculations we made using forbidden polynomials (cf.\ Lemma \ref{lem_taut}), we want to prove the following result about the kernel of the quotient map $j^*$.

\begin{thm} \label{thm_kernel}
	The kernel of $j^*: \HH \to \AA$ equals $\sum\limits_{p \geq 0,\ q > 0} \HH_p * (e_q^n \cup \HH_q)$.
\end{thm}

\begin{proof}
	Let $s_0: R_d \to \hat{R}_{d,n}$ be the zero section of $\hat{R}_{d,n}$ considered as a $G_d$-linear bundle on $R_d$. Under the identifications $\HH_d = A_{G_d}^*(R_d) \cong A_{G_d}^*(\hat{R}_{d,n}) \cong A_{G_d}^*(\pt)$, the map
	$$
		A_{G_d}^k(\pt) \cong A_{G_d}^k(R_d) \xto{}{(s_0)_*} A_{G_d}^{k+nd}(\hat{R}_{d,n}) \cong A_{G_d}^{k+nd}(\pt),
	$$
	is the multiplication with the top $G_d$-equivariant Chern class of $\hat{R}_{d,n}$. Identifying $A_{G_d}^*(\pt)$ with the ring of symmetric functions in $d$ variables, the top $G_d$-equivariant Chern class of $\smash{\hat{R}_{d,n}}$ is the $n$-th power of the $d$-th elementary symmetric polynomial. Taking this into account, the statement to prove is equivalent to showing that the horizontal sequence in the diagram
	\begin{center}
		\begin{tikzpicture}[description/.style={fill=white,inner sep=2pt}]
			\matrix(m)[matrix of math nodes, row sep=4.5em, column sep=2em, text height=1.5ex, text depth=0.25ex]
			{
				\bigoplus\limits_{\substack{p+q = d,\ q > 0\\ i+j = k - (m-1)pq}}\!\!\!\!\! \HH_p^i \otimes \HH_q^{j-nq} & \HH_d^k & \AA_d^k & 0 \\
				\bigoplus\limits_{\substack{p+q = d,\ q > 0\\ i+j = k - (m-1)pq}}\!\!\!\!\! \HH_p^i \otimes (e_q^n \cup \HH_q^j) & & & \\
			};
			\path[->, font=\scriptsize]
			(m-1-1) edge (m-1-2)
			(m-1-2) edge node[auto] {$j^*$} (m-1-3)
			(m-1-3) edge (m-1-4)
			($(m-1-1) + (1.5em,-2em)$) edge node[left] {$(s_0)_*$} ($(m-2-1) + (1.5em,1em)$)
			($(m-2-1) + (2.5em,1em)$) edge node[below] {$*$} (m-1-2)
			;
		\end{tikzpicture}
	\end{center}
	\vspace{1em}
	is exact. For all $p+q = d$, we have the K\"{u}nneth isomorphism
	$$
		\bigoplus_{i+j = k} \HH_p^i \otimes \HH_q^{j-nq} \xto{}{\cong} A_{G_p \times G_q}^{k-nq}(R_p \times R_q)
	$$
	and we have $ A_{G_p \times G_q}^{k-nq}(R_p \times R_q) \cong A_{G_p \times G_q}^{k-nq}(R_{\pqstar})$. Modulo these isomorphisms, writing $L_{p,q} = G_p \times G_q$ as we have already done before, we are interested in the map $\smash{A^{k-nq}_{L_{p,q}}(R_{\pqstar}) \to A^{k+(m-1)pq}_{G_d}(\hat{R}_{d,n})}$. These maps arise from the following morphisms
	\begin{center}
		\begin{tikzpicture}[description/.style={fill=white,inner sep=2pt}]
			\matrix(m)[matrix of math nodes, row sep=2em, column sep=1.5em, text height=1.5ex, text depth=0.25ex]
			{
				R_{\pqstar} \times^{L_{p,q}} U \\
				(M_{q \times n} \times R_{\pqstar}) \times^{L_{p,q}} U &[-.7em] &[-1em] \hat{R}_{\pqstar,n} \times^{P_{p,q}} U & \hat{R}_{d,n} \times^{P_{p,q}} U & \hat{R}_{d,n} \times^{G_d} U, \\
				\hat{R}_{\pqstar,n} \times^{L_{p,q}} U & & & \\
			};
			\path[->, font=\scriptsize]
			(m-1-1) edge node[auto] {$s_0$} (m-2-1)
			(m-3-1) edge (m-2-1)
			(m-3-1) edge (m-2-3)
			(m-2-3) edge node[auto] {$i$} (m-2-4)
			(m-2-4) edge node[auto] {$\pi$} (m-2-5)
			;
		\end{tikzpicture}
	\end{center}
	the non-horizontal ones being affine space bundles. Considering the cartesian squares
	\begin{center}
		\begin{tikzpicture}[description/.style={fill=white,inner sep=2pt}]
			\matrix(m)[matrix of math nodes, row sep=2em, column sep=4em, text height=1.5ex, text depth=0.25ex]
			{
				R_{\pqstar} \times^{L_{p,q}} U &	(M_{q \times n} \times R_{\pqstar}) \times^{L_{p,q}} U \\
				(M_{p \times n} \times R_{\pqstar}) \times^{L_{p,q}} U &	\hat{R}_{\pqstar,n} \times^{L_{p,q}} U \\
				(M_{p \times n} \times R_{\pqstar}) \times^{P_{p,q}} U &	\hat{R}_{\pqstar,n} \times^{P_{p,q}} U \\
			};
			\path[->, font=\scriptsize]
			(m-1-1) edge node[auto] {$s_0$} (m-1-2)
			(m-2-1) edge (m-1-1)
			(m-2-2) edge (m-1-2)
			(m-2-1) edge node[auto] {$s_0$} (m-2-2)
			(m-2-1) edge (m-3-1)
			(m-2-2) edge (m-3-2)
			(m-3-1) edge node[auto] {$s_0$} (m-3-2);
		\end{tikzpicture}
	\end{center}
	and using the commutativity of flat pull-back and proper push-forward, we are bound to show the exactness of
	\begin{equation}
		\bigoplus_{d = p+q,\ q > 0} A_{P_{p,q}}^{k-nq-(m-1)pq}(M_{p \times n} \times R_{\pqstar}) \to A_{G_d}^k(\hat{R}_{d,n}) \xto{}{j^*} A_{G_d}^k(\hat{R}_{d,n}^{\st}) \to 0. \tag{*}
	\end{equation}
	Let $(f,\phi)$ be an unstable point of $\hat{R}_{d,n}$. Then, the linear subspace $L(f,\phi)$ defined as the left-$\kk\langle\phi_1,\ldots,\phi_m \rangle$-submodule generated by $\im f$ is a proper subspace of $\kk^d$. Let $X_p$ be the $G_d$-invariant, closed subset of all $(f, \phi)$ where the subspace $L(f,\phi)$ has dimension at most $p$ (its natural scheme structure is the reduced one). This induces a filtration
	$$
		\hat{R}_{d,n}^{\mathrm{unst}} = X_{d-1} \supseteq X_{d-2} \supseteq \cdots \supseteq X_0 = R_d \times \{0\}
	$$
	by closed subsets. It corresponds to the Harder-Narasimhan stratification, as defined in \cite[Def.\ 3.3]{Reineke:03}. Denote $W_p := M_{p \times n} \times \smash{R_{\pqstar}}$ which we consider as a closed subset of $X_p$ by identifying $M_{p \times n}$ with the subspace of $M_{d \times n}$ of matrices whose image lies in the subspace of $\kk^d$ spanned by the first $p$ coordinate vectors. Then, $W_p$ carries a natural action of the parabolic $P_{p,q}$ (but it is not $G_d$-invariant). Evidently, the $G_d$-saturation of $W_p$ lies in $X_p$. Therefore, we obtain a morphism $\psi_p$ as the composition
	$$
		W_p \times^{P_{p,q}} U \to X_p \times^{P_{p,q}} U \to X_p \times^{G_d} U.
	$$
	As the first map is a closed immersion and the latter is a $G_d/P_{p,q}$-bundle, $\psi_p$ is proper. Consider the open subsets $X_p^o := X_p - X_{p-1}$ of $X_p$ and $W_p^o$ of $W_p$ defined as the subset of all $(f,\phi)$ such that $\dim L(f,\phi) = p$. We claim that $\psi_p$ induces an isomorphism
	$$
		W_p^o \times^{P_{p,q}} U \xto{}{\cong} X_p^o \times^{G_d} U.
	$$
	As $W_p^o \times^{P_{p,q}} U$ is naturally isomorphic to $(W_p^o \times^{P_{p,q}} G_d) \times^{G_d} U$ as a $G_d$-variety and as the map $\smash{X_p^o \times U \to X_p^o \times^{G_d} U}$ is a principal $G_d$-fiber bundle, it suffices by faithfully flat descent to show that
	$$
		W_p^o \times^{P_{p,q}} G_d \to X_p^o
	$$
	is an isomorphism of $G_d$-varieties. This will be done in Lemma \ref{lem_quot}. Denote by $W_p^c$ the complement of $W_p^o$ in $W_p$. Applying \cite[Ex.\ 1.8.1]{Fulton:98}, the cartesian diagram
	\begin{center}
		\begin{tikzpicture}[description/.style={fill=white,inner sep=2pt}]
			\matrix(m)[matrix of math nodes, row sep=2em, column sep=4em, text height=1.5ex, text depth=0.25ex]
			{
				W_p^c \times^{P_{p,q}} U &	W_p \times^{P_{p,q}} U \\
				X_{p-1} \times^{G_d} U &	X_p \times^{G_d} U \\
			};
			\path[->, font=\scriptsize]
			(m-1-1) edge (m-1-2)
			(m-1-1) edge node[auto] {$\pi'$} (m-2-1)
			(m-1-2) edge node[auto] {$\pi$} (m-2-2)
			(m-2-1) edge (m-2-2);
		\end{tikzpicture}
	\end{center}
	induces an exact sequence
	\begin{center}
		\begin{tikzpicture}[description/.style={fill=white,inner sep=2pt} scale=.5]
			\matrix(m)[matrix of math nodes, row sep=2em, column sep=1em, text height=1.5ex, text depth=0.25ex]
			{
				A_j(W_p^c \times^{P_{p,q}} U) &	A_j(X_{p-1} \times^{G_d} U) \oplus A_j(W_p \times^{P_{p,q}} U) &			A_j(X_p \times^{G_d} U) & 0 \\
				A_{r-pq}^{P_{p,q}}(W_p^c) &	 A_{r}^{G_d}(X_{p-1}) \oplus A_{r-pq}^{P_{p,q}}(W_p) &
									A_{r}^{G_d}(X_p) & 0,\\
			};
			\path[->, font=\scriptsize]
			(m-1-1) edge (m-1-2)
			(m-1-2) edge (m-1-3)
			(m-1-3) edge (m-1-4)
			(m-2-1) edge (m-2-2)
			(m-2-2) edge (m-2-3)
			(m-2-3) edge (m-2-4)
			;
			\draw[-] ($(m-1-1.south) + (-.1em,0)$) -- ($(m-2-1.north) + (-.1em,.2em)$);
			\draw[-] ($(m-1-1.south) + (.1em,0)$) -- ($(m-2-1.north) + (.1em,.2em)$);
			\draw[-] ($(m-1-2.south) + (-.1em,0)$) -- ($(m-2-2.north) + (-.1em,.2em)$);
			\draw[-] ($(m-1-2.south) + (.1em,0)$) -- ($(m-2-2.north) + (.1em,.2em)$);
			\draw[-] ($(m-1-3.south) + (-.1em,0)$) -- ($(m-2-3.north) + (-.1em,.2em)$);
			\draw[-] ($(m-1-3.south) + (.1em,0)$) -- ($(m-2-3.north) + (.1em,.2em)$);
			%
		\end{tikzpicture}
	\end{center}
	where the first map sends $\alpha$ to $\pi'_*\alpha + (- \alpha)$, the second maps $\beta + \beta'$ to $\beta + \pi_*\beta'$, and $r = j - \dim V + \dim G_d$. By induction on $p$, we obtain that the natural map
	$$
		A_r^{P_{0,d}}(W_0) \oplus \cdots \oplus A_{r - pq}^{P_{p,q}}(W_p) \to A_r^{G_d}(X_p)
	$$
	is onto. Inserting $p = d-1$ finally yields the exactness of the sequence
	\begin{center}
		\begin{tikzpicture}[description/.style={fill=white,inner sep=2pt}]
			\matrix(m)[matrix of math nodes, row sep=2em, column sep=1.5em, text height=1.5ex, text depth=0.25ex]
			{
				\bigoplus\limits_{d = p+q,\ q > 0} A^{P_{p,q}}_{r-pq}(W_p) &	A^{G_d}_r(\hat{R}_{d,n}) &	 A_r^{G_d}(\hat{R}_{d,n}^{\st}) &[-.8em] 0 \\
				\bigoplus\limits_{d = p+q,\ q > 0} A_{P_{p,q}}^{\dim W_p - (r-pq)}(W_p) &	A_{G_d}^{\dim \hat{R}_{d,n} -r}(\hat{R}_{d,n}) &	 A_{G_d}^{\dim \hat{R}_{d,n}^{\st} -r}(\hat{R}_{d,n}^{\st}) & 0. \\
			};
			\path[->, font=\scriptsize]
			(m-1-1) edge (m-1-2)
			(m-1-2) edge (m-1-3)
			(m-1-3) edge (m-1-4)
			(m-2-1) edge (m-2-2)
			(m-2-2) edge (m-2-3)
			(m-2-3) edge (m-2-4);
			\draw[-] ($(m-1-1.south) + (1.4em,-.2em)$) -- ($(m-2-1.north) + (1.4em,.4em)$);
			\draw[-] ($(m-1-1.south) + (1.6em,-.2em)$) -- ($(m-2-1.north) + (1.6em,.4em)$);
			\draw[-] ($(m-1-2.south) + (-.1em,0)$) -- ($(m-2-2.north) + (-.1em,.4em)$);
			\draw[-] ($(m-1-2.south) + (.1em,0)$) -- ($(m-2-2.north) + (.1em,.4em)$);
			\draw[-] ($(m-1-3.south) + (-.1em,0)$) -- ($(m-2-3.north) + (-.1em,.4em)$);
			\draw[-] ($(m-1-3.south) + (.1em,0)$) -- ($(m-2-3.north) + (.1em,.4em)$);
		\end{tikzpicture}
	\end{center}
	Choosing $r = \dim \hat{R}_{d,n} -k = md^2 +nd -k$ and observing that
	\begin{align*}
		\dim W_p - (r-pq) &= m(p^2+q^2 +pq) + np - (r - pq)\\
		&= md^2 + nd - nq - r - (m-1)pq \\
		&= k-nq-(m-1)pq,
	\end{align*}
	we have shown that the sequence (*) is exact.
\end{proof}

\begin{lem} \label{lem_quot}
	With the notation as in the proof of Theorem \ref{thm_kernel}, the natural map $W_p^o \times^{P_{p,q}} G_d \to X_p^o$ is an isomorphism.
\end{lem}

\begin{proof}
	Consider the morphism $L: X_p^o \to \Gr_{p,d} := \Gr_p(\kk^d) = G_d/P_{p,q}$ assigning to every point $(f,\phi) \in X_p^o$ the subspace $L(f,\phi)$. We show that
	\begin{center}
		\begin{tikzpicture}[description/.style={fill=white,inner sep=2pt}]
			\matrix(m)[matrix of math nodes, row sep=2em, column sep=4em, text height=1.5ex, text depth=0.25ex]
			{
				G_d \times W_p^o	&	X_p^o \\
				G_d		 	&	\Gr_{p,d} \\
			};
			\path[->, font=\scriptsize]
			(m-1-1) edge node[auto] {$\mathrm{act}$} (m-1-2)
			(m-1-1) edge node[auto] {$\pr$} (m-2-1)
			(m-1-2) edge node[auto] {$L$} (m-2-2)
			(m-2-1) edge (m-2-2);
		\end{tikzpicture}
	\end{center}
	is a cartesian diagram of varieties. In fact, $\smash{G_d \times_{\Gr_{p,d}} X_p^o}$ consists of those pairs $(g, (f,\phi))$ such that $L(v,\phi)$ equals the subspace of $\kk^d$ generated by $ge_1,\ldots,ge_p$. An isomorphism
	$$
		G_d \times W_p^o \to G_d \times_{\Gr_{p,d}} X_p^o
	$$
	is therefore given by mapping $(g, (f,\phi))$ to $(g, g \cdot (f,\phi))$.
\end{proof}

We deduce from Theorem \ref{thm_kernel} that the Chow ring of a non-commutative Hilbert scheme is tautologically presented. We make this statement a little more precise. For every $0 \leq p < d$, choose a basis $\BB_{p,q}$ of $\HH_p \otimes \HH_q$ as an $\HH_d$-module. It has cardinality $\binom{d}{p}$. Without loss of generality, we may assume that every basis element is a tensor product $f_{\lambda,p} \otimes g_{\lambda,q}$. Making this choice, we obtain:

\begin{cor}
	The kernel of $j^*: \HH_d = \Q[e_1,\ldots,e_d] \to A^*(\Hilb_{d,n})$, the homomorphism sending $e_\nu$ to the $\nu$-th Chern class of the universal bundle of $\Hilb_{d,n}$, is the ideal of $\HH_d$ generated by the expressions
	$$
		f_{\lambda,p} * (e_q^n \cup g_{\lambda,q})
	$$
	with $0 \leq p < d$, $q := d-p$ and $\lambda = 1,\ldots,\binom{d}{p}$.
\end{cor}

\begin{ex*}[continued]	
	Let's illustrate this result using once again our favorite non-commutative Hilbert scheme. Let $m=2$ and $d=3$ and $n=1$. We have $\HH_3 = \Q[e_1,e_2,e_3] = \Q[x,y,z]^{S_3}$.
	\begin{itemize}
		\item Let $p=0$. We obtain $\HH_p \otimes \HH_q = \HH_d$. Inserting $g = 1$ yields the relation $e_3 = 0$.
		\item For $p=1$, we get $\HH_p \otimes \HH_q = \Q[x][y,z]^{S_2}$. A basis as an $\HH_d$-module is given by $1,x,x^2$. Putting $f(x) = 1$ yields
		\begin{align*}
			0 &= yz(y-x)(z-x) + xz(x-y)(z-y) + xy(x-z)(y-z) \\
			&\equiv (xy + xz + yz)^2 \\
			&\equiv e_2^2
		\end{align*}
		when employing the relation $xyz = 0$. The other basis vectors result in multiples of $xyz$.
		\item Finally, let $p=2$. Then, a basis of $\HH_p \otimes \HH_q = \Q[x,y]^{S_2}[z]$ over $\HH_d$ is $1,z,z^2$. We consider $g(z) = 1$ first and obtain, using $xyz = 0$,
		$$
			0 \equiv e_1^3 - 4e_1e_2.
		$$
		After some lengthy computation, we see that for $g(z) = z$, we obtain the relation $e_1^4 = 0$. The basis element $z^2$ does not provide a new relation.
	\end{itemize}
	All in all, we have computed a presentation for the Chow ring of $\Hilb_{3,1}^{(2)}$. It is isomorphic to
	$$
		\Q[e_1,e_2]/(e_1^3 - 4e_1e_2, e_2^2, e_1^4).
	$$
\end{ex*}

\subsection{Two examples}

For $m=0$ and $m=1$, there is an explicit description of the CoHa as an exterior algebra and as a symmetric algebra, respectively, both over a vector space of countably infinite dimension (cf.\ \cite[Subsect.\ 2.5]{KS:11}). We would like to describe the ideal $\ker j^*$ under these isomorphisms.

Let $m=0$. In this case, the (non-commutative) Hilbert scheme $\Hilb_{d,n}$ is the Grassmannian $\Gr_{n-d}^n$ (i.e.\ empty for $d > n$). The multiplication in $\HH$ is given by the formula
\begin{align*}
	(f*g)(x_1,\ldots,x_d) = \sum_{\substack{1 \leq i_1 < \cdots < i_p \leq d\\1 \leq j_1 < \cdots < j_q \leq d\\\text{complementary}}} f(x_{i_1},\ldots,x_{i_p})g(x_{j_1},\ldots,x_{j_q}) \prod_{\mu = 1}^p \prod_{\nu = 1}^q (x_{j_\nu} - x_{i_\mu})^{-1}.
\end{align*}
It is easy to see that $(f * f)(x,y) = 0$ for all $f \in \HH_1$. Therefore, the embedding $\HH_1 \into \HH$ induces a homomorphism of (graded) algebras $\phi: \bigwedge^*(\HH_1) \to \HH$.
We identify the ring $\HH_1$ (equipped with the cup product $\cup$) with the polynomial ring $\Q[x]$. Let $\psi_0,\psi_1,\psi_2,\ldots$ be the basis of $\HH_1$ that corresponds to $x^0, x^1, x^2, \ldots$ under this isomorphism. Note that $\psi_i$ lives in bidegree $(1,2i+1)$. A basis of $\bigwedge^d(\HH_1)$ is given by expressions
$
	\psi_{k_1} \wedge \cdots \wedge \psi_{k_d},
$
where $k_1 < \cdots < k_d$ is an increasing sequence of $d$ non-negative integers. By induction on $d$, we can show that
$$
	(\psi_{k_1} * \cdots * \psi_{k_d})(x_1,\ldots,x_d) = s_\lambda(x_1,\ldots,x_d),
$$
where $s_\lambda$ is the Schur function belonging to the partition $\lambda = (k_d - d+1, \ldots, k_1)$. Hence, it follows that $\phi$ is an isomorphism. Let's determine $\II := \ker j^*$. Denoting $\II_d \sub \HH_d$ the $d$-th homogeneous component, Theorem \ref{thm_kernel} implies that $\II_d = \sum_{q=1}^d \HH_{d-q} * (e_q^n \cup \HH_q)$. We obtain that $\II_0 = 0$, and $\II_1 \sub \HH_1$ is $e_1^n \cup \HH_1$, which is generated by $\psi_n,\psi_{n+1},\ldots$ as a vector space.
%
In order to show that the surjection $\bigwedge^*(\psi_0,\ldots,\psi_{n-1}) \to \AA$ is an isomorphism, it suffices to prove that their generating series agree. We use Lemma \ref{lem_FJ} to compute the generating series of $\AA$. We have
\begin{align*}
	P_{\AA}(q,t) &:= \sum_d \sum_k (-1)^k \dim(\AA_{d,k}) q^{k/2}t^d \\
		&= \sum_d \sum_{k-d^2 \text{ even}} (-1)^k \dim\big( A^{(k-d^2)/2}(\Hilb_{d,n}^{(0)}) \big)q^{k/2}t^d \\
		&= \sum_d (-q^{1/2})^{d^2} \sum_j \dim\big( A^j(\Hilb_{d,n}^{(0)}) \big)q^j t^d \\
		&= \sum_d (-q^{1/2})^{d^2} \sum_{(j_1,\ldots,j_{n-d}) \in \smash{\mathcal{J}_{d,n}^{(0)}}} q^{\sum_\nu (d-j_\nu)} t^d.
\end{align*}
The set $\mathcal{J}_{d,n}^{(0)}$ is the set of non-decreasing sequences $0 \leq j_1 \leq \cdots \leq j_{n-d} \leq d$ (with no further restrictions). This also implies that $P_{\AA}$ is in fact a polynomial in $t$ of degree $n$. On the other hand, the generating series of the exterior algebra $\bigwedge^*(\psi_0,\ldots,\psi_{n-1})$ (with respect to its bigrading coming from the bidegrees of the $\psi_i$) is
\begin{align*}
	\prod_{i=0}^{n-1}(1-q^{i+1/2}t) &= \sum_{d=0}^n (-1)^d \sum_{1 \leq i_1 < \cdots < i_d \leq n} q^{\sum_l (i_l-1/2)}t^d.
\end{align*}
We have a bijection from $\mathcal{J}_{d,n}^{(0)}$ to the set of sequences $1 \leq i_1 < \cdots < i_d \leq n$ by assigning to $j_*$ the numbers $i_l = \sharp \{ \nu \mid j_\nu \leq d-l \} + l$. Under this bijection, we have that
$$
	\sum_l \left(i_l - \frac{1}{2}\right) = \sum_\nu (d-j_\nu) + (1+\cdots+d) - \frac{d}{2} = \sum_\nu (d-j_\nu) + \frac{d^2}{2}
$$
which proves that the two generating series coincide. We have proved:

\begin{cor}
	For $m = 0$, the $\HH = \bigwedge^*(\psi_0,\psi_1,\ldots)$-algebra $\AA = \bigoplus_d A^*(\Hilb_{d,n}^{(0)})$ equals the exterior algebra $\bigwedge\nolimits^{*}(\psi_0,\ldots,\psi_{n-1})$.
\end{cor}

Let's turn to the case $m=1$. 
The CoHa-multiplication has the form
$$
	(f*g)(x_1,\ldots,x_d) = \sum_{\substack{1 \leq i_1 < \cdots < i_p \leq d\\1 \leq j_1 < \cdots < j_q \leq d\\\text{complementary}}} f(x_{i_1},\ldots,x_{i_p})g(x_{j_1},\ldots,x_{j_q}).
$$
The immersion $\HH_1 \into \HH$ yields a homomorphism $\phi: \Sym^*(\HH_1) \to \HH$ of algebras (which respects the grading). This time, $\psi_i(x) = x^i$ in $\HH_1$ has bidegree $(1,2i)$. A basis element $\psi_{k_1}\cdots \psi_{k_d}$ of $\Sym^d(\HH_1)$ with $k_1 \geq \cdots \geq k_d$ is mapped to
$$
	(\psi_{k_1} * \cdots * \psi_{k_d})(x_1,\ldots,x_d) = c_\lambda m_\lambda(x_1,\ldots,x_d).
$$
In the above equation, $m_\lambda$ denotes the monomial symmetric function attached to the partition $\lambda = (k_1,\ldots,k_d)$ and $c_\lambda$ is some positive integer. This proves that $\phi$ is an isomorphism. We compute $\II$ and its inverse image under $\phi$. 
Like in the case $m = 0$, we see that $\psi_n,\psi_{n+1},\ldots$ span $\II_1$ as a vector space. We compare again the generating series of $\Q[\psi_0,\ldots,\psi_{n-1}]$ and $\AA$ to show that they are isomorphic. In analogy to the above computations, we find that
$$
	P_{\AA}(q,t) = \sum_d \sum_{(j_1,\ldots,j_n) \in \smash{\mathcal{J}_{d,n}^{(1)}}} q^{\sum_\nu (d-j_\nu)} t^d.
$$
The set $\mathcal{J}_{d,n}^{(1)}$ consists of all sequences $0 \leq j_1 \leq \cdots \leq j_n = d$. Moreover, the generating series of the symmetric algebra $\Q[\psi_0,\ldots,\psi_{n-1}]$ computes as 
$$
	\prod_{i=0}^{n-1} (1 - q^it)^{-1} = \sum_d \sum_{k_0 + \cdots + k_{n-1} = d} q^{\sum_i ik_i}t^d.
$$
A bijection between $\mathcal{J}_{d,n}^{(1)}$ and the set of all sequences $k_0,\ldots,k_{n-1}$ which sum to $d$ is given by assigning to $j_*$ the integers $k_{\nu-1} = j_\nu - j_{\nu-1}$. Then we see that
\begin{align*}
	\sum_{\nu=1}^n (\nu-1)k_{\nu-1} = (n-1)j_n - \sum_{\nu = 1}^{n-1} j_\nu = nd - \sum_{\nu=1}^n j_\nu = \sum_{\nu=1}^n (d-j_\nu)
\end{align*}
as $j_n = d$ by definition of $\mathcal{J}_{d,n}^{(1)}$. We obtain that the generating series agree and have thus proved the following:

\begin{cor}
	For $m=1$, the $\HH = \Q[\psi_0,\psi_1,\ldots]$-algebra $\AA = \bigoplus_d A^*(\Hilb_{d,n}^{(1)})$ coincides with the polynomial ring $\Q[\psi_0,\ldots,\psi_{n-1}]$.
\end{cor}

	%
	\bibliographystyle{abbrv}
	\bibliography{Literature}
\end{document}